\documentclass[11pt,a4paper]{amsart}
\usepackage{amsfonts, amsmath}
\usepackage[pdftex]{graphics}
\usepackage{psfrag}
\usepackage{tikz}

\newtheorem{theorem}{Theorem}[section]

\newtheorem{lemma}[theorem]{Lemma}
\newtheorem{cor}[theorem]{Corollary}

\numberwithin{equation}{section}

\def\proof{\noindent {\bf Proof. }}
\def\eop{\hfill {\hfill $\Box$} \medskip}

\newcommand{\cov}{\mathrm{cov}}
\newcommand{\Ker}{\mathrm{Ker}}
\newcommand{\val}{\mathrm{val}}

\newcommand{\R}{{\mathcal{R}}}
\newcommand{\cH}{{\mathcal{H}}}
\newcommand{\M}{{\mathcal{M}}}
 \newcommand{\Y}{{\mathcal{Y}}}
\newcommand{\rank}{\operatorname{rank}}

\newcommand{\sm}{\setminus}

\begin{document}

\title{\bf Rigid cylindrical frameworks with two
coincident points}
\author[B. Jackson]{B. Jackson}
\address{School of Mathematical Sciences\\ Queen Mary, University of London\\
E1 4NS \\ U.K. } \email{b.jackson@qmul.ac.uk}
\author[V. E. Kaszanitzky]{V. E. Kaszanitzky}
\address{Department of Mathematics and Statistics\\ Lancaster University\\
LA1 4YF \\ U.K. }
\email{viktoria@cs.elte.hu}
\author[A. Nixon]{A. Nixon}
\address{Department of Mathematics and Statistics\\ Lancaster University\\
LA1 4YF \\ U.K. }
\email{a.nixon@lancaster.ac.uk}
\date{\today}

\begin{abstract}
We develop a rigidity theory for frameworks in $\mathbb{R}^3$ which
have two coincident points but are otherwise generic and only
infinitesimal motions which are tangential to a family of cylinders
induced by the realisation are considered.
We then apply our results to show that vertex splitting, under the additional assumption that the new edge is redundant, preserves the property of being generically globally
rigid on families of concentric cylinders.
\end{abstract}

\maketitle



\section{Introduction}

A \emph{framework} $(G,p)$ in $\mathbb{R}^d$ is the combination of a finite, simple graph $G=(V,E)$ and a map $p:V\rightarrow \mathbb{R}^d$. It is \emph{rigid} if every edge-length preserving continuous motion of the vertices arises as a congruence of $\mathbb{R}^d$ (see, for example, \cite{Whlong} for basic definitions and background). The theory of generic rigidity aims to characterise the graphs $G$ for which $(G,p)$ is rigid for all generic choices of  $p$.  This was accomplished  by Laman \cite{laman} for $d=2$, but is a long-standing open problem for $d\geq 3$.

We are interested in frameworks in $\mathbb{R}^3$ whose motions are restricted such that all vertices are realised on a fixed surface and only those continuous motions that keep $(G,p)$ on the surface are considered. Generic rigidity in this context has been characterised for various surfaces \cite{NOP,NOP2}. In this paper we consider frameworks on concentric cylinders in which two of the vertices are mapped to the same point, but are otherwise generic. For such frameworks we give a combinatorial characterisation of rigidity.

Frameworks in $\mathbb{R}^2$ with two coincident points were studied in \cite{FJK} where the following deletion-contraction characterisation of rigidity was proved. A graph \(G\) is \(uv\)-rigid in $\mathbb{R}^2$ if
there exists a realisation $p$ of $G$ in $\mathbb{R}^2$ such that $p(u)=p(v)$, $p|_{V-v}$ is generic and the framework $(G,p)$ is rigid.
We use \(G-uv\) to denote the graph formed from \(G\) by deleting the edge \(uv\) if it exists and \(G/uv\) to denote the graph which arises from \(G\) by contracting the vertices \(u\) and \(v\).

\begin{theorem}\cite{FJK} \label{thm:planeuv}
Let \(G=(V,E)\) be a graph and let \(u,v\in V\) be distinct vertices. Then \(G\) is \(uv\)-rigid in $\mathbb{R}^2$ if and only if \(G-uv\) and \(G/uv\) are both rigid in $\mathbb{R}^2$.
\end{theorem}

Our main result is an analogue of this theorem for frameworks on concentric cylinders. One motivation for studying frameworks on cylinders with coincident points is an ongoing work to understand global rigidity for frameworks on concentric cylinders \cite{JMN,JNstress}. In particular we utilise our main result to prove that vertex splitting, under the additional assumption that the new edge is redundant, preserves generic global rigidity for families of concentric cylinders.

In Section \ref{sec:surfaces} we provide background on frameworks on concentric cylinders. In Section \ref{countsec} we define a count matroid $\mathcal{M}_{uv}(G)$ on a graph $G$ with two distinguished vertices $u$ and $v$. In Section \ref{hennsec} we derive an inductive construction for graphs whose edge set is independent in $\mathcal{M}_{uv}(G)$. We then use this construction to prove our characterisation of rigidity on families of concentric cylinders for frameworks in which $u$ and $v$ are coincident but are otherwise generic. In Section \ref{sec:global} we discuss global rigidity and apply our coincident point result to prove that the vertex splitting operation preserves generic global rigidity for frameworks on families of concentric cylinders when the new edge is redundant. Finally, in Section \ref{sec:further} we comment on extensions to other surfaces.

In this paper we consider simple graphs only as parallel edges correspond to the same distance constraint and thus one of them is always redundant.

\section{Frameworks on concentric cylinders}
\label{sec:surfaces}

Let $G=(V,E)$ where $V=\{v_1,\dots,v_n\}$. We will consider realisations of $G$ on a family of concentric cylinders $\Y=\mathcal{Y}_1 \cup \mathcal{Y}_2 \cup \dots \cup \mathcal{Y}_k$ where $\mathcal{Y}_i=\{(x,y,z)\in \mathbb{R}^3:x^2+y^2=r_i\}$ and $r=(r_1,\dots,r_k)$ is a vector of positive real numbers.
A \emph{framework} $(G,p)$ on $\mathcal{Y}$ is an ordered pair consisting of a graph $G$ and a realisation $p$ such that $p(v_i)\in \mathcal{Y}$ for all $v_i\in V$.

Two frameworks $(G,p)$ and $(G,q)$ on $\mathcal{Y}$ are \emph{equivalent} if $\|p(v_i)-p(v_j)\|=\|q(v_i)-q(v_j)\|$ for all edges $v_iv_j\in E$. Moreover $(G,p)$ and $(G,q)$ on $\mathcal{Y}$ are \emph{congruent} if $\|p(v_i)-p(v_j)\|=\|q(v_i)-q(v_j)\|$ for all pairs of vertices $v_i,v_j\in V$.
The framework $(G,p)$ is \emph{rigid} on $\mathcal{Y}$ if there exists an $\epsilon>0$ such that every framework $(G,q)$ on $\mathcal{Y}$ which is equivalent to $(G,p)$, and has $\| p(v_i)-q(v_i)\|<\epsilon$ for all $1\leq i \leq n$, is congruent to $(G,p)$. Moreover $(G,p)$ is \emph{minimally rigid} if $(G,p)$ is rigid but $(G-e,p)$ is not for any $e\in E$.
The framework $(G,p)$ is \emph{generic} on $\mathcal{Y}$ if td$[\mathbb{Q}(r,p):\mathbb{Q}(r)]=2n$.

It was shown in \cite{NOP} that a generic framework $(G,p)$ on any family of concentric cylinders is rigid if and only if it is infinitesimally rigid in the following sense. An \emph{infinitesimal flex} $s$ of $(G,p)$ on $\mathcal{Y}$ is a map $s:V\rightarrow \mathbb{R}^3$ such that $s(v_i)$ is tangential to $\mathcal{Y}$ at $p(v_i)$ for all $v_i\in V$ and $(p(v_j)-p(v_i))\cdot (s(v_j)-s(v_i))=0$ for all $v_jv_i\in E$. The framework $(G,p)$ is \emph{infinitesimally rigid} on $\mathcal{Y}$ if every infinitesimal flex is an infinitesimal isometry of $\mathbb{R}^3$.

The \emph{rigidity matrix} $R_{\Y}(G,p)$ is the $(|E|+|V|)\times 3|V|$ matrix
\[ R_{\Y}(G,p)=\begin{pmatrix}R_3(G,p)\\ S(G,p) \end{pmatrix}\]
where:
$R_3(G,p)$ has rows indexed by $E$ and 3-tuples of columns indexed by $V$ in which, for $e=v_iv_j\in E$, the submatrices in row $e$ and columns $v_i$ and $v_j$ are $p(v_i)-p(v_j)$ and $p(v_j)-p(v_i)$, respectively, and all other entries are zero; $S(G,p)$ has rows indexed by $V$ and 3-tuples of columns indexed by $V$ in which, for $v_i\in V$, the submatrix in row $v_i$ and column $v_i$ is $\bar p(v_i)=(x_i,y_i,0)$ when $p(v_i)=(x_i,y_i,z_i)$.
The \emph{rigidity matroid} $\mathcal{R}^{\mathcal{Y}}(G)$ is the row matroid of $R_{\Y}(G,p)$ for any generic $p$.

A graph $G=(V,E)$ is \emph{$(2,2)$-sparse} if $|E'|\leq 2|V'|-2$ for all subgraphs $(V',E')$ of $G$. Moreover $G$ is \emph{$(2,2)$-tight} if $G$ is $(2,2)$-sparse and $|E|=2|V|-2$.

The following characterisation of generic rigidity on $\mathcal{Y}$ was proved in \cite{NOP}.

\begin{theorem}\label{thm:cylinderlaman}
Let $(G,p)$ be a generic framework on a union of concentric cylinders $\mathcal{Y}$. Then $(G,p)$ is minimally rigid if and only if $G$ is a complete graph on at most 3 vertices or $G$ is $(2,2)$-tight and simple.
\end{theorem}


\subsection{Coincident realisations on concentric cylinders}

Let $G=(V,E)$ be a graph and suppose $u,v\in V$. A framework $(G,p)$ on $\Y$ is \emph{$uv$-coincident} if $p(u)=p(v)$. A \emph{generic $uv$-coincident framework} is a $uv$-coincident framework $(G,p)$ for which $(G-u,p|_{V-u})$ is generic. We denote the $uv$-coincident cylinder rigidity matroid by $\mathcal{R}_{uv}^{\Y}(G)$ (this is the row matroid of $R_\Y(G,p)$ for any generic $uv$-coincident realisation $(G,p)$).
Note that the matroid depends on $G$ but not on the choice of generic $uv$-coincident realisation. That is, for any two generic $uv$-coincident realisations $(G,p)$ and $(G,p')$ on $\Y$, we get the same matroid.
We also use $r_{uv}(G)$ to denote the rank of $\mathcal{R}_{uv}^{\Y}(G)$.
We say that $G$ is \emph{$uv$-rigid} if  $r_{uv}(G)=2|V|-2$ and that $G$ is \emph{minimally $uv$-rigid} if $G$ is $uv$-rigid and $|E|=2|V|-2$.

Note that the term $uv$-rigid and the notation $r_{uv}(G)$ refer to generic realisations on a fixed family of concentric cylinders $\Y$, and hence appear to depend on $\Y$. We will see, however, that this is not the case since our characterisation of $\mathcal{R}^\Y_{uv}(G)$ depends only on the graph $G$.

\section{A count matroid}
\label{countsec}

In this section we define a count matroid $\M_{uv}(G)$ on the edge
set of a graph $G$ with two distinguished vertices $u$ and $v$. We
will show that  $\M_{uv}(G)$ is equal to $\mathcal{R}^\Y_{uv}(G)$ in
Section \ref{hennsec}.

Let $G=(V,E)$ be a graph.  For some $X\subseteq V$ let $G[X]$ denote
the subgraph of $G$ induced by $X$ and let $E_G(X)$ be the set of
edges of $G[X]$. Thus $i_G(X)=|E_G(X)|$. For a family
\(\mathcal{S}=\{S_{1},S_{2},\dots ,S_{k}\}\), where \(S_{i}\subseteq
V\) for all \(i=1,\dots,k\), we define
\(E_G(\mathcal{S})=\cup_{i=1}^{k}E_G(S_{i})\) and put
\(i_G(\mathcal{S})=|E_G(\mathcal{S})|\). We also define
$\cov({\mathcal{S}})=\{ (x,y) : x,y\in V, \{x,y\}\subseteq S_i\
\hbox{for some}\ 1\leq i\leq k \}.$ We say that ${\mathcal{S}}$ {\it
covers} a set $F\subseteq E$ if $F\subseteq \cov({\mathcal{S}})$.
The degree of a vertex $w$ is denoted by $d_G(w)$.
We let $N_G(w)=\{ z\in V : wz\in E \}$ denote the {\it neighbours} of $w$ in $G$. We may omit the subscripts referring to $G$ if the graph is clear from the context.

Let \(G=(V,E)\) be a graph and $u,v\in V$ be  two distinct vertices
of $G$. Let $\mathcal{H}=\{H_1,...,H_k\}$ be a family with
$H_i\subseteq V$, $1\leq i\leq k$. We say that $\mathcal{H}$ is
\emph{\(uv\)-compatible} if \(u,v\in H_{i}\) and \(|H_{i}|\geq 3\)
hold for all $1\leq i\leq k$. See Figure \ref{fig:ex} for an example. We define the {\it value} of subsets
of $V$ and of \(uv\)-compatible families as follows. For a nonempty
subset $H\subseteq V$,
we let
$$\val(H)=2|H|-t_H,$$ where $t_H=4$ if $H= \{u,v\}$, $t_H=3$ if $H\neq \{u,v\}$ and $|H| \in \{2,3\}$,
and $t_H=2$ otherwise.
We will often denote \(t_{H_i}\) by \(t_i\) for short. For a
\(uv\)-compatible family \(\mathcal{H}=\{H_{1},H_{2},\dots
,H_{k}\}\) we let
$$\val(\mathcal{H})=\sum_{i=1}^{k}\val(H_{i})-2(k-1).$$ Note that if
\(\mathcal{H}=\{H\}\) is a \(uv\)-compatible family containing only
one set then the two definitions are compatible, i.e.
$\val(\mathcal{H})=\val(H)$ holds.

We  say that $G$ is {\it $uv$-sparse} if for all $H\subseteq V$ with
$|H|\geq 2$ we have $i_{G}(H)\leq \val(H)$ and for all
\(uv\)-compatible families $\mathcal{H}$ we have
$i_{G}(\mathcal{H})\leq \val(\mathcal{H})$. Note that if $G$ is
$uv$-sparse then \(G\) is simple and $uv\notin E$ must hold. A set \(H\subseteq V\) of
vertices with \(|H|\geq 2\) (resp. a \(uv\)-compatible family
\(\mathcal{H}=\{H_{1},\dots,H_{k}\}\)) is called \emph{tight} if
\(i_G(H)=\val(H)\) (resp. $i_G(\mathcal{H})=\val(\mathcal{H})$)
holds. We will show that the edge sets of the $uv$-sparse subgraphs
of $G$ form the independent sets of a matroid that we will denote by $\M_{uv}(G)$.

The next lemmas will enable us to 'uncross' tight sets and tight $uv$-compatible families in a sparse graph. The first result follows immediately from the definition of the $i$- and $\val$- functions.

\begin{lemma} \label{lem:2sets}
Let \(X,Y\subseteq V\) be distinct vertex sets in \(G\). Then\\
(a) $i(X)+i(Y)\leq i(X\cup Y)+i(X\cap Y)$ and\\
(b) if \(X\cap Y\neq\emptyset\),
then
$ \val(X)+\val(Y) +t_{X} +t_{Y} = \val(X\cup Y)+ \val(X\cap Y)+t_{X\cup Y} +t_{X\cap Y}.$
\end {lemma}

\begin{lemma}\label{lem:intersect3}
Let \(\mathcal{H}=\{H_{1},\dots,H_{k}\}\) be a \(uv\)-compatible family in $G$.\\
(a) Suppose \(|H_{i}\cap H_{j}|\geq 3\) for some pair \(1\leq
i<j\leq k\). Then  there is a \(uv\)-compatible family
\(\mathcal{H}'\)
with $\cov(\mathcal{H})\subseteq \cov(\mathcal{H}')$ and
\(\val(\mathcal{H}')<\val(\mathcal{H})\). \\
(b) Suppose $G$ is $uv$-sparse and $\mathcal{H}$ is tight. Then $H_i\cap H_j=\{u,v\}$ for all $1\leq i \leq k$.
\end{lemma}

\proof
(a) We may assume that \(i=k-1\), \(j=k\).
Let \(\mathcal{H}'=\{H_{1},\dots,H_{k-2},H_{k-1}\cup H_{k}\}\). Using Lemma \ref{lem:2sets}(b) we have
$\val(H_{k-1})+\val(H_k) \geq\val( H_{k-1} \cup H_{k}) + \val( H_{k-1} \cap  H_{k}) $. Hence
\begin{align*} \val(\mathcal{H})&=& \sum_{l=1}^{k}\val(H_{l})-2(k-1)=\sum_{l=1}^{k-2}\val(H_{l})-2((k-1)-1)+\val(H_{k-1})+\val(H_{k})-2\\
&\geq&\sum_{l=1}^{k-2}\val(H_{l})+\val(H_{k-1}\cup H_{k})-2((k-1)-1)+\val (H_{k-1}\cap H_k) -2>\val(\mathcal{H}').\end{align*}
Clearly, we have $\cov(\mathcal{H})\subseteq \cov(\mathcal{H}')$.
\\
(b) Since  $\mathcal{H}$ is tight, if \(|H_{i}\cap H_{j}|\geq 3\)
for some pair \(1\leq i<j\leq k\) then, by (a), we have
$\val(\mathcal{H}')<\val(\mathcal{H})=i(\mathcal{H})\leq
i(\mathcal{H}')$. This contradicts the $uv$-sparsity of $G$. Hence
$H_i\cap H_j=\{u,v\}$ for all $1\leq i \leq k$. \eop

\begin{lemma}\label{lem:2size4}
Let \(\mathcal{H}=\{H_{1},\dots,H_{k}\}\) be a \(uv\)-compatible family with $H_i\cap H_j=\{u, v\}$ for all $1\leq i<j \leq k$ and \(|H_{k}|\geq 4\). Then \(\mathcal{H}'=\{H_{1},\dots,H_{k-2},H_{k-1}\cup H_{k}\}\)
 is a \(uv\)-compatible family with $\cov(\mathcal{H})\subset\cov(\mathcal{H}')$ and for which \(\val(\mathcal{H}')\leq\val(\mathcal{H})+1\) with equality only if $|H_{k-1}|=3$.
Furthermore,  if $G$ is $uv$-sparse, $\mathcal{H}$ is tight and $|H_{k-1}|\geq 4$, then $\mathcal{H}'$ is tight.
\end{lemma}

\proof Using Lemma \ref{lem:2sets}(b) and the facts that \(t_k=t_{H_{k-1}\cup H_k}=2\) and \(t_{H_{k-1}\cap  H_{k}}=4\) we have
$\val(H_{k-1})+\val(H_k)=\val( H_{k-1} \cup H_{k}) + \val( H_{k-1} \cap  H_{k})+4-t_{k-1}= \val( H_{k-1} \cup H_{k}) +4-t_{k-1}$.
Hence
\begin{eqnarray*}
\val(\mathcal{H})=\sum_{l=1}^{k}\val(H_{l})-2(k-1)&=&\sum_{l=1}^{k-2}\val(H_{l})-2((k-1)-1)+\val(H_{k-1})+\val(H_{k})-2\\
&=&  \sum_{l=1}^{k-2}\val(H_{l})+\val(H_{k-1}\cup H_{k})-2((k-1)-1)+2-t_{k-1}\\
&=&\val(\mathcal{H}')+2-t_{k-1}.
\end{eqnarray*}
Thus \(\val(\mathcal{H}')\leq\val(\mathcal{H})+1\) with equality only if $|H_{k-1}|=3$. Clearly, we have $\cov(\mathcal{H})\subset\cov(\mathcal{H}')$.

Now suppose $G$ is $uv$-sparse, $\mathcal{H}$ is tight and $|H_{k-1}|\geq 4$. Then $\val(\mathcal{H}')\leq \val(\mathcal{H})=i(\mathcal{H})=i(\mathcal{H}')$, so $\mathcal{H}'$ is tight.
\eop

\begin{lemma} \label{lem:sets}
Let $G=(V,E)$ be $uv$-sparse and let \(X,Y\subseteq V\) be tight
sets in \(G\) with \(X\cap Y\neq\emptyset\) and  $|X|,|Y|\geq4$.
Then $|X\cap Y|\not\in \{2,3\}$ and \(X\cup Y\) and \(X\cap Y\) are
both tight.
\end{lemma}
\proof
We have
\[2|X|-2+2|Y|-2=i(X)+i(Y)\leq i(X\cup Y)+i(X\cap Y)\]
\[\leq 2|X\cup Y|-t_{X\cup Y}+2|X\cap Y|-t_{X\cap Y}=2|X|+2|Y|-2-t_{X\cap Y}.\]
This implies that \(t_{X\cap Y}=2\) and equality holds throughout.
Thus \(X\cup Y\) and \(X\cap Y\) are both tight and either \(|X\cap
Y|\geq4\) or \(|X\cap Y|=1\). \eop

\begin{lemma}\label{lem:set2intersection}
Let \(\mathcal{H}=\{H_{1},\dots,H_{k}\}\)  be a \(uv\)-compatible
family with $H_j\cap H_l=\{u,v\}$ for all $1\leq j<l\leq k$, and let
\(Y\subseteq V\) be a set of vertices with \(|Y|\geq4\), and $|Y\cap
\{u,v\}|\leq 1$. Suppose that for some \(1\leq i\leq
k\) either  \(|Y\cap H_{i}|\geq 2\), or  \(|Y\cap H_{i}|=1\) and $|H_i|\geq 4$.
Then there is a \(uv\)-compatible family \(\mathcal{H}'\) with
\(\cov(\mathcal{H})\cup \cov(Y)\subseteq\cov(\mathcal{H}')\) and
\(\val(\mathcal{H}')\leq \val(\mathcal{H})+\val(Y)\).
Furthermore, if \(G\) is \(uv\)-sparse and \(\mathcal{H}\) and \(Y\)
are both tight then \(\mathcal{H}'\) and \(Y\cap H_i\) are also
tight.
\end{lemma}

\proof
Let $S=\{H_i\in \cH\,:\, |Y\cap H_{i}|\geq 2 \mbox{ or } |Y\cap H_{i}|=1 \mbox{ and } |H_i|\geq 4\}$.
Renumbering the sets of $\mathcal{H}$, if necessary, we may assume
that $S=\{H_i\in \cH\,:\,j\leq i\leq k$\}, for some \(j\leq k\). Let
\(X=Y\cup(\cup_{i=j}^{k}H_i)\) and
\(\mathcal{H}'=\{H_{1},\dots,H_{j-1},X\}\). Then
\(\cov(\mathcal{H})\cup \cov(Y)\subseteq \cov(\mathcal{H}')\) and
\[|X|=\sum_{i=j}^{k}|H_{i}|+|Y|-2(k-j)-\sum_{i=j}^{k}|H_{i}\cap Y|+|Y\cap \{u,v\}|(k-j).\]
This gives
\begin{eqnarray*}\val(\mathcal{H})+\val(Y)&=&\sum_{i=1}^{k}\val(H_{i})-2(k-1)+\val(Y)\\
&=&\sum_{i=1}^{j-1}\val(H_{i})-2(j-1)+\sum_{i=j}^{k}(2|H_{i}|-t_i)-2(k-j)+(2|Y|-2)\\
&=&\sum_{i=1}^{j-1}\val(H_{i})+(2|X|-2)-2(j-1)+4(k-j)-\sum_{i=j}^{k}t_{H_i}\\
&&+2\sum_{i=j}^{k}|Y\cap H_{i}|-2(k-j)-2|Y\cap \{u,v\}|(k-j)\\
&\geq&\sum_{i=1}^{j-1}\val(H_{i})+\val(X)-2(j-1)+
\sum_{i=j}^{k}(2|Y\cap H_{i}|-t_{H_i}).
\end{eqnarray*}
If $|Y\cap H_{i}|\geq 2$ then $\val(Y\cap H_{i})=2|Y\cap H_{i}|-t_{Y\cap H_i}\leq 2|Y\cap H_{i}|-t_{H_i}$. On the other hand,
if $|Y\cap H_{i}|=1 \mbox{ and } |H_i|\geq 4$, then $t_{Y\cap H_i}=2=t_{H_i}$ and we have $\val(Y\cap H_{i})= 2|Y\cap H_{i}|-t_{H_i}$. Thus, in both cases,
\[\val(\mathcal{H})+\val(Y)\geq \val(\mathcal{H}')+\sum_{i=j}^{k} \val(Y\cap H_i)\]
and so \(\val(\mathcal{H}')\leq \val(\mathcal{H})+\val(Y)\).

Now, suppose that $G$ is $uv$-sparse and \(\mathcal{H}\) and \(Y\) are tight. Then we have
\[i(\mathcal{H}')+\sum_{i=j}^{k}i(Y\cap H_i)\geq i(\mathcal{H})+i(Y)=
\val(\mathcal{H})+ \val(Y)\geq \]
\[ \geq \val(\mathcal{H}')+\sum_{i=j}^{k} \val(Y\cap H_i)\geq i(\mathcal{H}')+\sum_{i=j}^{k}i(Y\cap H_i),\]
where the first inequality follows from the fact that edges spanned by \(\mathcal{H}\) or \(Y\) are spanned by \(\mathcal{H}'\) and if some edge is spanned by both \(\mathcal{H}\) and \(Y\) then it is spanned by \(Y\cap H_i\) for some \(i\). The equality holds because \(\mathcal{H}\) and \(Y\) are tight, and the second inequality holds by our calculations above. The last inequality holds because \(G\) is \(uv\)-sparse. Hence equality must hold everywhere, which implies
that \(\mathcal{H}'\) is tight and that \(Y\cap H_i\) is also tight for all $j\leq i\leq k$.
\eop

\begin{lemma}\label{lem:Y2}
Let \(\mathcal{H}=\{H_{1},\dots,H_{k}\}\) be a \(uv\)-compatible family with $H_i\cap H_j=\{u,v\}$ for all $1\leq i< j\leq k$, and let $Y\subseteq V$ be a set of vertices with \(|Y|\geq4\), $Y\cap \{u,v\}=\emptyset$ and $|Y\cap H_i|\leq 1$ for all $1\leq i\leq k$. Suppose that \(|Y\cap H_{i}|=|Y\cap H_{j}|=1\) for some pair \(1\leq i<j\leq k\). Then there is a \(uv\)-compatible family \(\mathcal{H}'\) with \(\cov(\mathcal{H})\cup \cov(Y)\subseteq \cov(\mathcal{H}')\) for which \(\val(\mathcal{H}')\leq \val(\mathcal{H})+\val(Y)\).
Furthermore, if \(G\) is \(uv\)-sparse and \(\mathcal{H}\) and \(Y\) are both tight, then $\mathcal{H}'$ is tight and $|H_{i}|=|H_j|=3$.
\end{lemma}
\proof
We may assume that \(i=k-1\) and \(j=k\). Let
\(\mathcal{H}'=\{H_{1},\dots,H_{k-2},H_{k-1}\cup H_{k}\cup Y\}\). We have \(\cov(\mathcal{H})\cup \cov(Y)\subseteq \cov(\mathcal{H}')\) and
\begin{eqnarray*}
\val(\mathcal{H})+\val(Y)&=&\sum_{i=1}^{k}\val(H_{i})-2(k-1)+\val(Y)\\
&=&\sum_{i=1}^{k-2}\val(H_{i})-2((k-1)-1)-2+\val(H_{k-1})+\val(H_{k})+\val(Y).
\end{eqnarray*}
Using Lemma \ref{lem:2sets}(b) twice and the fact that $|H_{k-1}\cap (H_k\cup Y)|=3$ we obtain
\begin{eqnarray*}
\val(H_{k-1})+\val(H_{k})+\val(Y)&=&\val(H_{k-1})+\val(H_{k}\cup Y)+2-t_{H_k}\\
&=&\val(H_{k-1}\cup H_{k}\cup Y)+8-t_{H_{k-1}}-t_{H_k}\\
&\geq&\val(H_{k-1}\cup H_{k}\cup Y)+2,
\end{eqnarray*}
with equality only if $|H_{k-1}|=|H_k|=3$. Thus \(\val(\mathcal{H}')\leq \val(\mathcal{H})+\val(Y)\) as claimed.

Now suppose that \(G\) is \(uv\)-sparse.
and \(\mathcal{H}\) and \(Y\) are both tight.
Then we have
\[i(\mathcal{H})+i(Y)=
\val(\mathcal{H})+ \val(Y)\geq \val(\mathcal{H}')\geq i(\mathcal{H}')\geq i(\mathcal{H})+i(Y)\]
where the last inequality follows since
\(|Y\cap H_{k-1}|=|Y\cap H_{k}|=1\) and
$|Y\cap H_i|\leq 1$ for all $1\leq i\leq k$. Hence equality must hold throughout. Thus \(\mathcal{H}'\)  is tight and $|H_{k-1}|=|H_k|=3$.
\eop

\begin{lemma} \label{tightsystem}
Let \(G=(V,E)\) be $uv$-sparse and suppose that there is a tight \(uv\)-compatible family in \(G\). Then there is a unique tight \(uv\)-compatible family
\(\mathcal{H}_{\max}\) in \(G\) for which $\cov(\mathcal{H})\subseteq \cov(\mathcal{H}_{\max})$ for all tight \(uv\)-compatible families \(\mathcal{H}\) of \(G\). In addition, if  $\cH_{\max}=\{X_1,X_2,\ldots,X_k\}$ and $|X_1|\geq |X_2|\geq\ldots\geq |X_k|$, then:\\
(a) $X_i\cap X_j=\{u,v\}$ for all $1\leq i<j\leq k$;\\
(b) $|X_i|=3$ for all $2\leq i\leq k$;\\
(c) $N(u,v)\subseteq V(\cH_{max})$.\\
Furthermore, if $Y\subseteq V$ is tight, $|Y|\geq 4$, $\cov(Y)\not \subseteq \cov(\cH_{\max})$, and $Y\cap X_i\neq \emptyset$ for some $1\leq i\leq k$, then
$|Y\cap X_i|=1$, $|X_i|=3$, $Y\cap \{u,v\}=\emptyset$,  and $Y\cap X_j=\emptyset$ for all $j\neq i$.
\end{lemma}

\proof
Let \(\mathcal{H}_1=\{X_{1},X_{2},\dots,X_{k}\}\) be a tight \(uv\)-compatible family in $G$ labeled such that $|X_1|\geq |X_2|\geq\ldots\geq |X_k|$ and suppose that $\cov(\cH_1)$ is maximal with respect to inclusion.
Then Lemmas \ref{lem:intersect3} and \ref{lem:2size4} imply that  \(X_{i}\cap X_{j}=\{u,v\}\) holds for all $1\leq i<j\leq k$ and \(|X_i|=3\) for all \(2\leq i\leq k\). Suppose for a contradiction that \(\mathcal{H}_{2}=\{Y_{1},Y_{2},\dots,Y_{l}\}\) is another tight \(uv\)-compatible family
whose cover is maximal, labeled so that $|Y_1|\geq |Y_2|\geq\ldots\geq |Y_{l}|$.
We will use notation \(X_i=\{u,v,x_i\}\) for \(2\leq i\leq k\) and \(Y_j=\{u,v,y_i\}\) for \(2\leq j\leq l\). Without loss of generality we can assume that if \(|X_1|=|Y_1|=3\) then \(X_1\neq Y_1\).

We define two \(uv\)-compatible families as follows: let
\[\mathcal{H}_{\cap}=\{Z\subseteq V:|Z|\geq3 \mbox{ and $X_{i}\cap Y_{j}=Z$ for some $X_i\in \mathcal{H}_1$, $Y_j\in \mathcal{H}_2$}\};\]
let
\[\mathcal{H}_{\cup} =\{X_1\cup Y_1\}\cup\{X_i:2\leq i\leq k \mbox{ and }x_i\not\in X_1\cup Y_1\}\cup\{Y_j:2\leq j\leq l\mbox{ and }y_j\not\in X_1\cup Y_1\}\]
if \(|X_1\cap Y_1|\geq3\), and
\[\mathcal{H}_{\cup}=\{X_1\}\cup\{Y_1\}\cup\{X_i:2\leq i\leq k\mbox{ and }x_i\not\in X_1\cup Y_1\}\cup\{Y_j:2\leq j\leq l\mbox{ and }y_j\not\in X_1\cup Y_1\}\]
if \(|X_1\cap Y_1|=2\).

It is easy to see that $\mathcal{H}_{\cup}$ and $\mathcal{H}_{\cap}$ are both $uv$-compatible. For convenience we rename the families as
\(\mathcal{H}_{\cup}=\{A_{1},\dots,A_{p}\}\) and
\(\mathcal{H}_{\cap}=\{B_{1},\dots,B_{q}\}\), where $A_1=X_1\cup Y_1$ and $B_1=X_1\cap Y_1$ if \(|X_1\cap Y_1|\geq3\), and $A_1=X_1$ and $A_2=Y_1$ if \(|X_1\cap Y_1|=2\). It follows from their construction that \(|A_i|=3\) for all \(3\leq i\leq p\) and \(|B_j|=3\) for all \(2\leq j\leq q\) and also at least one of \(|A_2|=3\), \(|B_1|=3\) holds. It can be seen easily that \(p+q=k+l\). We also have \(i(\mathcal{H}_{1})+i(\mathcal{H}_{2})\leq i(\mathcal{H}_{\cup})+i(\mathcal{H}_{\cap})\), since the family
\(\mathcal{H}_{\cup}\) spans all the edges spanned by \(\mathcal{H}_{1}\)
or \(\mathcal{H}_{2}\) and \(\mathcal{H}_{\cap}\) spans all the edges spanned by both \(\mathcal{H}_{1}\) and \(\mathcal{H}_{2}\). Thus
\begin{eqnarray*}
\val(X_1)+3(k-1)-2(k-1)+\val(Y_1)+3(l-1)-2(l-1)=
\val(\mathcal{H}_{1})+\val(\mathcal{H}_{2})\\
=i(\mathcal{H}_{1})+i(\mathcal{H}_{2})\leq
i(\mathcal{H}_{\cup})+i(\mathcal{H}_{\cap})\leq
\val(\mathcal{H}_{\cup})+
\val(\mathcal{H}_{\cap})\\
=\val(A_1)+\max\{\val(A_2),\val(B_1)\}+3(p-1)-2(p-1)+3(q-1)-2(q-1).
\end{eqnarray*}

We will show that equality occurs at both ends of the above inequality. Since \(k-1+l-1=p-1+q-1\), it will suffice to show that \(\val(X_1)+\val(Y_1)\geq\val(A_1)+\max\{\val(A_2),\val(B_1)\}\). This is immediate if \(|X_1\cap Y_1|=2\) and follows from Lemma \ref{lem:2sets}(b) when \(|X_1\cap Y_1|\geq 3\).

Hence equality must hold throughout the displayed inequality. In particular,  \(\mathcal{H}_{\cup}\) and
\(\mathcal{H}_{\cap}\) are both tight. Since
\(\cov(\mathcal{H}_{1})\cup
\cov(\mathcal{H}_{2})\subseteq\cov(\mathcal{H}_{\cup})\), the
maximality of the covers implies that  \(\cov(\mathcal{H}_{1})=
\cov(\mathcal{H}_{2})\) which in turn gives
\(\mathcal{H}_{1}=\mathcal{H}_{2}\).

We have now shown that $\mathcal{H}_{1}=\mathcal{H}_{\max}$ is unique and that properties (a) and (b) hold. To see that (c) holds choose $x\in N(u,v)$ and suppose that $x\not\in V(\cH_{max})$. Let $\cH'=\cH_{\max}+\{u,v,x\}$.  Then $i(\cH')\geq i(\cH_{\max})+1$ and
$\val(\cH')=\val(\cH_{\max})+1$, so $\cH'$ is tight and hence contradicts the maximality of $\cH_{\max}$.

To complete the proof we suppose that $Y\subseteq V$ is tight,
$|Y|\geq 4$, $\cov(Y)\not \subseteq \cov(\cH_{\max})$, and $Y\cap
X_i\neq \emptyset$ for some $1\leq i\leq k$. If $\{u,v\}\subseteq Y$
then $\cH=\{Y\}$ would be a $uv$-compatible family with
$\cov(\cH)\not\subseteq \cov(\cH_{\max})$. This would contradict the
maximality of $\cH_{\max}$ and hence $\{u,v\}\not\subseteq Y$. If
$|Y\cap X_i|\geq 2$ or $|Y\cap X_i|=1$ and $|X_i|\geq 4$ then Lemma
\ref{lem:set2intersection} would imply that there exists a
$uv$-compatible family $\cH'$  with $\cov(\cH_{\max})\cup
\cov(Y)\subseteq \cov(\cH')$. Hence $|Y\cap X_i|\leq 1$ and
$|X_i|=3$. This tells us that $| Y\cap X_j| \leq 1 $ for all $ j $
and hence $\cov (Y)\cap \cov (\cH_{\max})=\emptyset $. If $Y\cap
\{u,v\}\neq \emptyset$ then putting $\cH'=\cH_{\max}\cup \{Y\cup
\{u,v\}\}$ we have $i(\cH')\geq i(\cH)+2|Y|-2$ and
$\val(\cH')=\val(\cH)+2|Y|-2$, so $\cH'$ would contradict the
maximality of $\cH_{\max}$. Thus $Y\cap \{u,v\}= \emptyset$. If
$Y\cap X_j\neq \emptyset$ for some $j\neq i$ then Lemma \ref{lem:Y2}
now gives us a tight $uv$-compatible family $\cH'$ with
$\cov(\cH_{\max})\cup \cov(Y)\subseteq \cov(\cH')$. Hence $Y\cap
X_j=\emptyset$ for all $j\neq i$.
 \eop

Note that Lemma \ref{tightsystem} tells us in particular that if $G$ is $uv$-sparse and $Y\subseteq V$ is tight with $\{u,v\}\cap Y\neq \emptyset$, then $Y\subseteq X_i$ for some $X_i\in \mathcal{H}_{\max}$.

\subsection{The matroid and its rank function}

We first remind the reader of the simple $(2,2)$-sparse matroid and its rank function. Given a graph $G=(V,E)$, a set $F\subseteq E$ is independent if and only if it is simple and induces a $(2,2)$-sparse subgraph. A system $\mathcal{K}=\{ H_{1},\dots,H_{k} \}$ of subsets of $V$ is \emph{thin} if \(|H_{i}\cap H_{j}|\leq 1\) for all pairs \(1\leq i,j\leq k\) with equality only if $|H_i|=2$ or $|H_j|=2$. The value of the system $\mathcal{K}$ is given by $\sum_{H_i\in\mathcal{\mathcal{K}}}\val(H_i)$.

Now we define the count matroid $\mathcal{M}_{uv}(G)$.
Let $G=(V,E)$ be a graph and $u,v\in V$ be distinct vertices of $G$. We will prove that the family of sets
\begin{equation}
\label{I}
\mathcal{I}_G=\{F: F\subseteq E \mbox{ and } (V,F)\ \hbox{is}\ uv\hbox{-sparse} \}
\end{equation}
defines a matroid  $\mathcal{M}_{uv}(G)$ on $E$ and characterise the rank function of this matroid. We need the following definition.

Let $\mathcal{H}=\{X_{1},\dots,X_{t}\}$ be a $uv$-compatible family and let $H_{1},\dots,H_{k}$ be subsets of $V$ of size at least two.
The system $\mathcal{K}=\{ \mathcal{H}, H_{1},\dots,H_{k} \}$ is a \emph{$uv$-cover} of $F\subseteq E$ if $F\subseteq \cov(\mathcal{H}) \cup \cov(\{H_1,\dots,H_k\})$. It is \emph{thin} if \\
(i) $\{H_{1},\dots,H_{k}\}$ is thin,\\
(ii) \(X_{i}\cap X_{j}=\{u,v\}\) for all pairs \(1\leq i,j\leq t\), and\\
(iii) \(|H_{i}\cap X_{j}|\leq 1\) for all \(1\leq i\leq k\), \(1\leq j\leq t\).\\
The value of the system $\mathcal{K}$ is given by \(\val(\mathcal{K})=\val (\mathcal{H})+ \sum_{i=1}^{k}\val(H_{i})\).

We will show that the rank of an arbitrary subset $F\subseteq E$ in $\M_{uv}(G)$ is given by
\begin{equation}
\label{rank}
r(F)=\min\{ \val (\mathcal{K}):\mathcal{K}\hbox{ is a thin cover of }F\}.
\end{equation}

 \begin{figure}[htp]
\begin{center}
   \begin{tikzpicture}[very thick,scale=.45]
\filldraw (-4,0) circle (3pt)node[anchor=east]{$u$};
\filldraw (4,0) circle (3pt)node[anchor=west]{$v$};
\filldraw (0,4) circle (3pt)node[anchor=south]{$v_1$};
\filldraw (0,2.5) circle (3pt)node[anchor=north]{$v_2$};
\filldraw (0,-2.5) circle (3pt)node[anchor=south]{$v_3$};
\filldraw (-1.5,-4) circle (3pt)node[anchor=north]{$v_4$};
\filldraw (1.5,-4) circle (3pt)node[anchor=north]{$v_5$};

\draw[black,thick]
(1.5,-4) -- (-1.5,-4) -- (0,-2.5) -- (1.5,-4) -- (4,0) -- (0,4) -- (-4,0) -- (0,2.5) -- (0,4);

\draw[black,thick]
(0,2.5) -- (4,0) -- (0,-2.5) -- (-4,0) -- (-1.5,-4);

\end{tikzpicture}
\end{center}
\vspace{-0.3cm}
\caption{An example of a simple $(2,2)$-tight graph \(G\) which is not independent in $\mathcal{M}_{uv}(G)$.}
\label{fig:ex}
\end{figure}
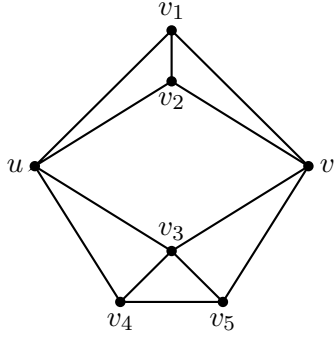

Let \(G=(V,E)\) be the graph shown in Figure \ref{fig:ex}. It is not difficult to see that \(G\) is (2,2)-sparse and simple, and hence $E$ is independent in the simple (2,2)-sparse matroid. We will show that \(E\) is not independent in $\mathcal{M}_{uv}(G)$. Consider the following sets: \(X_1=\{u,v,v_1\}\), \(X_2=\{u,v,v_2\}\) and \(X_3=\{u,v,v_3,v_4,v_5\}\). Then \(\mathcal{H}=\{X_1,X_2,X_3\}\) is a \(uv\)-compatible family of \(G\) with \(\val(\mathcal{H})=\val(X_1)+\val(X_2)+\val(X_3)-2\cdot2=(2\cdot3-3)+(2\cdot3-3)+(2\cdot5-2)-4=10\) and \(\cov(\mathcal{H})=E-v_1v_2\). Hence
$i_G(\mathcal{H})=11>\val(\mathcal{H})$
so \(E\) is dependent in $\M_{uv}(G)$.

\begin{theorem} \label{matroid}
Let $G=(V,E)$ be a graph and $u,v\in V$ be distinct vertices of $G$. Then $\mathcal{M}_{uv}(G)=(E,\mathcal{I}_G)$ is a matroid on ground-set $E$, where $\mathcal{I}_G$ is defined by (\ref{I}). The rank of a set \(E'\subseteq E\) in $\mathcal{M}_{uv}(G)$  is equal to
$$\min \{ \val(\mathcal{K}): \mathcal{K}\
\hbox{ is either a thin cover or a thin $uv$-cover of } E'\}.$$
\end{theorem}

\proof Let $\mathcal{I}=\mathcal{I}_G$, let $E'\subseteq E$ and let \(F\subseteq E'\) be a maximal subset of \(E'\) in \(\mathcal{I}\). Since \(F\in \mathcal{I}\) we have \(|F|\leq \val(\mathcal{K})\) for all ($uv$-)covers \(\mathcal{K}\) of \(E'\). We shall prove that there is a thin ($uv$-)cover \(\mathcal{K}\) of \(E'\) with \(|F|=\val(\mathcal{K})\), from which the theorem will follow.

Let \(J=(V,F)\) denote the subgraph defined by the edge set \(F\). First suppose that there is no tight \(uv\)-compatible family in \(J\) and consider the following cover of $F$:
\[\mathcal{K}_{1}=\{H_{1},H_{2},\dots,H_{k}\},\]
where \(H_{1},H_{2},\dots,H_{t}\) are the maximal tight sets with size at least four in \(J\) for some \(t\leq k\) and \(H_{t+1},\dots,H_{k}\) are the pairs of end vertices of edges in $J'=(V,F-\cup_{i=1}^{t}E(H_i))$. Clearly $\mathcal{K}_1$ is a cover of $F$. It is thin by Lemma \ref{lem:sets}. Thus
\[|F|=\sum_{j=1}^{k}|E_{J}(H_{j})|=\sum_{j=1}^{k}(2|H_{j}|-t_j)=\val(\mathcal{K}_{1})\]
follows. We claim that \(\mathcal{K}_{1}\) is a cover of \(E'\). To see this consider an edge \(ab=e\in E'-F\). Since \(F\) is a maximal subset of \(E'\) in \(\mathcal{I}\) we have \(F+e\not\in \mathcal{I}\). By our assumption there is no tight \(uv\)-compatible family in \(J\), and hence there must be a tight set \(X\) in \(J\) with  \(a,b\in X\). Hence \(X\subseteq H_{i}\) for some \(1\leq i\leq k\) which implies that  $\mathcal{K}_1$ covers $e$.

Next suppose that there is a tight \(uv\)-compatible family in $J$ and consider the following $uv$-cover of \(F\):
\[\mathcal{K}_{2}=\{\mathcal{H}_{\max},H_{1},H_{2},\dots,H_{k}\},\]
where \(\mathcal{H}_{\max}=\{X_{1},X_{2},\dots,X_{l}\}\) is the \(uv\)-compatible family of $G$ for which $\cov(\mathcal{H}_{\max})$ is maximal (given by Lemma \ref{tightsystem}) and \(H_{1},H_{2},\dots,H_{t}\) are the maximal tight sets with size at least four of $J'=(V,F-E(\mathcal{H}_{\max}))$ for some \(t\leq k\) and \(H_{t+1},\dots,H_{k}\) are the pairs of end vertices of edges in $J''=(V,F-E(\mathcal{H}_{\max})-\cup_{i=1}^{t}E(H_i))$. Then $\mathcal{K}_{2}$ is a $uv$-cover of $F$. By Lemmas  \ref{lem:sets} and \ref{tightsystem}, the $uv$-cover \(\mathcal{K}_{2}\) is thin, and hence
\[|F|=\sum_{i=1}^{l}|E_{J}(X_{i})|+\sum_{j=1}^{k}|E_{J}(H_{j})|=\sum_{i=1}^{l}(2|X_{i}|-t_i)-2(l-1)+\sum_{j=1}^{k}(2|H_{j}|-t_j)
=\val(\mathcal{K}_{2}).\]
We claim that \(\mathcal{K}_{2}\) is a $uv$-cover of \(E'\). As above, let \(ab=e\in E'-F\) be an edge. By the maximality of $F$ we have \(F+e\not\in \mathcal{I}\). Thus either there is a tight set \(X\subseteq V\) in $J$ with \(a,b\in \cov(X)\) or there is a tight \(uv\)-compatible family \(\mathcal{H}'=\{Y_{1},\dots,Y_{t}\}\) in $J$ with \(a,b\in Y_{i}\) for some \(1\leq i\leq t\).

In the latter case Lemma \ref{tightsystem} implies that \(\cov(\mathcal{H}')\subseteq \cov(\mathcal{H}_{\max})\) and hence $e$ is covered by \(\mathcal{K}_{2}\).
In the former case, when \(a,b\in X\) for some tight set \(X\) in $J$, we have $|X|\geq 5$ since if $|X|=2,3$ or $4$ then $X$ induces a complete graph in $J$ and $e=ab$ would be an edge of $F$. Lemma \ref{tightsystem} now gives \(|X\cap \cup_{i=1}^{l} X_{i}|\leq 1\). Then $E(X)\subseteq E(J')$ and hence \(X\subseteq H_{i}\) for some \(1\leq i\leq k\), since every edge of $J'$ induces a tight set and every tight set is contained in a maximal tight set. Thus $e$ is covered by \(\mathcal{K}_{2}\), as claimed.
\eop

\section{Characterisation of the $uv$-coincident cylinder rigidity matroid}
\label{hennsec}
Our aim is to show that the $uv$-coincident cylinder rigidity matroid $\R_{uv}^\Y(G)$ of a graph $G=(V,E)$ is equal to the count matroid $\M_{uv}(G)$.
To simplify terminology we will say that $G$ {\em is independent in  $\R_{uv}^\Y$, respectively $\M_{uv}$}, if $E$ is independent in $\R_{uv}^\Y(G)$, respectively $\M_{uv}(G)$.

We first show that independence in $\R_{uv}^\Y$ implies independence in $\M_{uv}$. Let $G/uv$ denote the graph obtained from \(G\) by contracting the vertex pair $u,v$ into a new vertex $z_{uv}$ (and deleting the resulting loops and parallel copies of edges).
Given a $uv$-coincident realisation $(G,p)$ of $G$ on $\Y$ we obtain a realisation $(G/uv,p_{uv})$ of $G/uv$ on $\Y$ by putting $p_{uv}(z_{uv})=p(u)=p(v)$ and $p_{uv}(x)=p(x)$ for all $x\in V-\{u,v\}$. Furthermore, each vector in the kernel of $R_{\Y}(G/uv,p_{uv})$ determines a vector in the kernel of $R_{\Y}(G,p)$ in a natural way. It follows that

\begin{equation}
\label{eq:kernel}
\dim \Ker R_{\Y}(G,p) \geq\dim\Ker R_{\Y}(G/uv,p_{uv}).
\end{equation}

We can use this fact to prove that independence in $\R_{uv}^\Y$ implies independence in \(\mathcal{M}_{uv}\).

\begin{lemma}
\label{lem:necessary}
Let \(G=(V,E)\) be a graph and let $u,v\in V$ be distinct vertices. If $G$ is independent in \(\mathcal{R}_{uv}^\Y\) then $G$ is independent in \(\mathcal{M}_{uv}\).
\end{lemma}

\proof Let $(G,p)$ be an independent $uv$-coincident realisation of $G$. Independence implies that \(i(X)\leq \val(X)\) holds for all \(X\subseteq V\) by Theorem \ref{thm:cylinderlaman}. Since $p(u)=p(v)$, $uv\notin E$ follows.

Let \(\mathcal{H}=\{X_{1},\dots,X_{k}\}\) be a \(uv\)-compatible family and consider the subgraph \(F=(\cup_{i=1}^{k}X_{i},\cup_{i=1}^{k}E(X_{i}))\). By contracting the vertex pair $u,v$ in $F$ we obtain the graph $F/uv$, in which \(\mathcal{H}_{uv}=\{X_{1}/uv,\dots,X_{k}/uv\}\) is a cover where \(X_{i}/uv\) denotes the set that we get from \(X_{i}\) by identifying \(u\) and \(v\). Thus we get \(r(F/uv)\leq\sum_{i=1}^{k}\val(X_{i}/uv)=\sum_{i=1}^{k}(2(|X_{i}|-1)-t_i)\) by using (\ref{rank}). This bound and (\ref{eq:kernel}) imply that $\dim\Ker R_{\Y}(F,p)\geq\dim\Ker R_{\Y}(F/uv,p_{uv})\geq 2(|\cup_{i=1}^{k}X_{i}|-1)-\sum_{i=1}^{k}(2|X_{i}|-(t_i+2))$. Since $(G,p)$ is an independent $uv$-coincident realisation of $G$, we have
$$i_F({\mathcal{H}})=|F|\leq
2\left|\bigcup_{i=1}^{k}X_{i}\right|-\left(2\left(|\bigcup_{i=1}^{k}X_{i}|-1\right)-\sum_{i=1}^{k}\left(2|X_{i}|-(t_i+2)\right)\right)=$$
$$ \sum_{i=1}^{k}(2|X_{i}|-t_i)-2(k-1)=\val({\mathcal{H}}).$$
Thus \(G\) is independent in \(\mathcal{M}_{uv}\), as claimed. \eop

We next define operations on $uv$-sparse graphs and use them to show that independence in $\M_{uv}$ implies independence in $\R_{uv}^\Y$.

The (two-dimensional versions of) the well-known Henneberg operations are as follows. Let $G=(V,E)$ be a graph. The {\it $0$-extension} operation (on a pair of distinct vertices $a,b\in V$) adds a new vertex $z$ and two edges $za,zb$ to $G$. The {\it 1-extension} operation (on edge $ab\in E$ and vertex $c\in V-\{a,b\}$)
deletes the edge $ab$, adds a new vertex $z$ and edges $za,zb,zc$.

We shall need the following specialized versions. Let $u,v\in V$ be two distinct vertices. The \emph{0-\(uv\)-extension} operation is a $0$-extension on a pair $a,b$ with \(\{a,b\}\neq \{u,v\}\). The \emph{1-\(uv\)-extension} operation is a $1$-extension on some edge $ab$ and vertex $c$ for which $\{u,v\}$ is not a subset of $\{a,b,c\}$. The inverse operations are called \emph{0-\(uv\)-reduction} and \emph{1-\(uv\)-reduction}, respectively.

We will also need two further moves. The \emph{vertex-to-$K_4$} move deletes a vertex $w$ and substitutes in a copy of $K_4$ with $V(K_4)\cap V(G)=\{w\}$ and with an arbitrary replacement of edges $xw$ by edges $xy$ with $y\in V(K_4)$. The inverse operation is known as a \emph{$K_4$-contraction}. A \emph{vertex-to-4-cycle} move takes a vertex $w$ with neighbours $v_1,v_2,\dots,v_k$ for any $k\geq 2$, splits $w$ into two new vertices $w,w'$ with $w'\notin V(G)$, adds edges $wv_1,w'v_1,wv_2,w'v_2$ and then arbitrarily replaces edges $xw$ with edges $xy$ where $x\in \{v_3,\dots,v_k\}$ and $y\in \{w,w'\}$. The inverse move is known as a
\emph{4-cycle-contraction}. The only difference in the specialised versions of these moves are that we require $|V(K_4)\cap \{u,v\}|\leq 1$ in a $uv$-$K_4$-contraction and similarly $|V(C_4)\cap \{u,v\}|\leq 1$ in a $uv$-4-cycle-contraction.

We first consider the 0-extension and 1-extension operations. It was shown in \cite{NOP} that these operations preserve independence in $\mathcal{R}^\Y$. The same arguments can be used to verify analogous results for $\mathcal{R}^\Y_{uv}$.

\begin{lemma} \label{lem:ext}
Let \(G=(V,E)\) be independent in $\R_{uv}^\Y$ and suppose that $G'$ is obtained from $G$ by a 0-\(uv\)-extension or a 1-\(uv\)-extension Then \(G'\) is independent in $\R_{uv}^\Y$.
\end{lemma}

In the case of 0-extensions we will also need the following result.

\begin{lemma}\label{lem:0ext}
Let $(G,p)$ be a generic realisation of a graph $G=(V,E)$ and $v\in V$. Suppose that $R_\Y(G,p)$ has linearly independent rows. Let $G'$ be obtained by performing a 0-extension which adds a new vertex $u$ to $G$. Put \(p'(a)=p(a)\) for all \(a\in V\), and put \(p'(u)=p(v)\). Then $R_{\Y}(G',p')$ has linearly independent rows.
\end{lemma}

\proof
The 0-extension adds 3 rows and 3 columns to
$R_\Y(G,p)$, the 3 columns being 0 everywhere except the 3 new rows.
The genericness of $p$ and the fact that $uv\notin E$ implies the
new $3\times 3$ block is invertible. Hence $R_\Y(G',p')$ has linearly
independent rows so \(G'\) is independent in $\mathcal{R}_{uv}^\Y$.
\eop

We next consider the vertex-to-4-cycle operation. It was shown in \cite{NOP2} that this operation preserves  independence in $\mathcal{R}^{\Y}$. A similar argument would yield the analogous result for $\mathcal{R}^{\Y}_{uv}$ but we will need a stronger result that a vertex-to-4-cycle move which creates two coincident vertices preserves independence in $\mathcal{R}^{\Y}$.

\begin{lemma}\label{lem:4cycleind}
Suppose $(G,p)$ is a framework on $\Y$, $R_\Y(G,p)$ has linearly independent rows and $w\in V$ with neighbours  $v_1,v_2,\dots,v_k$. Suppose further that $p(w)-p(v_1), p(w)-p(v_2)$ and $\bar p(w)$ are linearly independent where $\bar p(w)$ is the projection of $p(w)$ onto the plane $z=0$. Let $G'$ be obtained by performing a vertex-to-4-cycle operation at $w$ in $G$ such that $v_1$ and $v_2$ are both adjacent to $w$ and $w'$ in $G'$. Put \(p'(a)=p(a)\) for all \(a\in V-w\) and put \(p'(w)=p'(w')=p(w)\). Then $R_{\Y}(G',p')$ has linearly independent rows.
\end{lemma}

\proof
We will construct $R_{\Y}(G',p')$ from $R_{\Y}(G,p)$ by a series of
simple matrix operations that preserve the independence of the rows.

We first add three zero columns corresponding to \(w'\).
We then add three rows corresponding to the edges \(w'v_1,w'v_2\)
and the vertex \(w'\). Adding these rows increases the rank by 3
since $p(w)-p(v_1), p(w)-p(v_2)$ and $\bar p(w)$ are linearly
independent so the \(3\times3\) matrix formed by the entries in the
columns corresponding to \(w'\) and the rows corresponding to
\(w'v_1,w'v_2,w'\) is non-singular and the rest of the entries in
these columns are zero. The matrix \(M\) we obtain by this
modification has the following form: 
\[\begin{tabular}{c|ccccccc|}
  \multicolumn{1}{ c }{}  & \multicolumn{3}{ c }{\(\overbrace{\hspace{13mm}}^{w}\)}  &
  \multicolumn{3}{ c }{\(\overbrace{\hspace{13mm}}^{w'}\)} & \multicolumn{1}{ c }{} \\
  \cline{2-8}
  \((wv_1)\)  & \multicolumn{3}{ c }{\(p(w)-p(v_1)\)}  & \multicolumn{3}{ c }{\(\bf 0\)}& \(\star\) \\
  \((wv_2)\)  & \multicolumn{3}{ c }{\(p(w)-p(v_2)\)}  & \multicolumn{3}{ c }{\(\bf 0\)}& \(\star\) \\
& \multicolumn{3}{ c }{\(\vdots\)}  & \multicolumn{3}{ c }{\(\vdots\)}& \(\vdots\)\\
  \((wv_i)\)  & \multicolumn{3}{ c }{\(p(w)-p(v_i)\)}  & \multicolumn{3}{ c }{\(\bf 0\)}& \(\star\) \\
& \multicolumn{3}{ c }{\(\vdots\)}  & \multicolumn{3}{ c }{\(\vdots\)}& \(\vdots\)\\
  \((w'v_1)\)  &  \multicolumn{3}{ c }{\(\bf 0\)} & \multicolumn{3}{ c }{\(p(w)-p(v_1)\)} & \(\star\) \\
  \((w'v_2)\)  &  \multicolumn{3}{ c }{\(\bf 0\)}  & \multicolumn{3}{ c }{\(p(w)-p(v_2)\)} & \(\star\) \\
& \multicolumn{3}{ c }{\(\vdots\)}  & \multicolumn{3}{ c }{\(\vdots\)}& \(\vdots\)\\
  \cline{2-8}
  \(w\)  & \multicolumn{3}{ c }{\(\bar p(w)\)} &\multicolumn{3}{ c }{\(\bf 0\)} & \(\bf 0\) \\
  \(w'\)  & \multicolumn{3}{ c }{\(\bf 0\)} & \multicolumn{3}{ c }{\(\bar p(w)\)}& \(\bf 0\)\\
& \multicolumn{3}{ c }{\(\vdots\)}  & \multicolumn{3}{ c }{\(\vdots\)}& \(\vdots\)\\
  \cline{2-8}
\end{tabular}=M\]

To obtain $R_{\Y}(G',p')$ from $M$ we need  to modify some of the rows in \(M\) corresponding to edges \((wv_i)\) into the form of rows corresponding to edges \((w'v_i)\), i.e. we need to move the entries in the columns of \(w\) to the columns of \(w'\) and replace them with zeros. We will do this one by one.

Since \((p(w)-p(v_1))\), \((p(w)-p(v_2))\) and \(\bar p(w)\) are linearly independent, for every \(3\leq i\leq k\) there exist unique values \(\alpha,\beta,\gamma\) such that \(\alpha(p(w)-p(v_1))+\beta(p(w)-p(v_2))+\gamma\bar p(w)=(p(w)-p(v_i))\). Now subtract the row of \((wv_1)\) multiplied by \(\alpha\), the row of \((wv_2)\) multiplied by \(\beta\) and the row of \(w\) multiplied by \(\gamma\) from the row of \((wv_i)\) in \(M\). Then add the row of \((w'v_1)\) multiplied by \(\alpha\), the row of \((w'v_2)\) multiplied by \(\beta\) and the row of \(w'\) multiplied by \(\gamma\) to the same row (and change its label from
$(wv_i)$ to $(w'v_i)$) for every neighbour \(v_i\) of \(w'\) in \(G'\) to obtain $R_{\Y}(G',p')$. These operations also preserve independence, thus we conclude that the rows of $R_{\Y}(G',p')$ are independent.
\eop

\begin{cor} \label{cor:4cycleuv}
Let \(G\) be  independent in $\R_{uv}^\Y$ and suppose
that $G'$ is obtained from $G$ by a vertex-to-4-cycle operation.
Then \(G'\) is  independent in $\R_{uv}^\Y$.
\end{cor}

\proof
We choose a generic $uv$-coincident realisation $(G,p)$. Then $(G,p)$ satisfies the hypotheses of Lemma \ref{lem:4cycleind}. Hence $G'$ has a $uv$-coincident realisation $(G',p')$ such that $R_{\Y}(G',p')$ has linearly independent rows. It follows that every generic $uv$-coincident realisation is independent.
\eop

We next consider a generalisation of the vertex-to-$K_4$ operation. It was shown in \cite{NOP} that this operation preserves independence in $\mathcal{R}^{\Y}$. We will need an analogous result for $uv$-coincident realisations.

\begin{lemma}\label{lem:extension}
Let $G=(V,E)$ be a graph with $|E|=2|V|-2$ and let $u,v\in V$ be distinct vertices. Suppose $H\subset G$ is chosen so that either:\\
(a) $u,v\in V(H)$, $H$ is minimally $uv$-rigid on $\Y$ and $G/H$ is minimally rigid on $\Y$, or\\
(b) $|\{u,v\}\cap V(H)|\leq 1$, $H$ is minimally rigid on $\Y$ and $G/H$ is minimally $uv$-rigid on $\Y$. (Taking $z$ to be the vertex of $G/H$ obtained by contracting $H$ when $\{u,v\}\cap V(H)=z$.)  \\
Then $G$ is $uv$-rigid on $\Y$.
\end{lemma}

\proof
(a) Let $|V|=n$, $|V(H)|=r$ and consider $R_\Y(G,p)$ where $(G,p)$ is a generic $uv$-coincident framework on $\Y$ and $p=(p(v_1),p(v_2),\dots,p(v_n))$. By reordering rows and columns if necessary we can write $R_\Y(G,p)$ in the form
\[\begin{pmatrix} R_\Y(H,p|_H) & 0 \\ M_1(p) & M_2(p) \end{pmatrix} \]
where $M_2(p)$ is a square matrix with $3(n-r)$ rows.

Suppose, for a contradiction, that $G$ is not $uv$-rigid. Then there exists a vector $m\in \ker R_\Y(G,p)$ which is not an infinitesimal isometry of $\Y$. Since $(H,p|_H)$ is $uv$-rigid we may suppose that $m=(0,\dots,0,m_{r+1},\dots,m_n)$. Consider the realisation $(G,p')$ where $p'=(p(v_r), p(v_r), \dots,p(v_r), p(v_{r+1}), \dots, p(v_n))$ and define the realisation $(G/H,p^*)$ by setting $p^*=(p(v_r),p(v_{r+1}),\dots,p(v_n))$.  Since $p^*$ is generic, $(G/H,p^*)$ is infinitesimally rigid on $\Y$ by assumption.

Now, $M_2(p)$ is square with the nonzero vector $(m_{r+1},\dots,m_n)\in \ker M_2(p)$. Hence $\rank M_2(p) < 3(n-r)$. Since $p$ is generic, we also have $\rank M_2(p') < 3(n-r)$ and hence there exists a nonzero vector $m'\in \ker M_2(p')$. Therefore we have
\[ \begin{pmatrix} R_\Y(G/H,p^*)\end{pmatrix} \begin{pmatrix}0 \\ m' \end{pmatrix}= \begin{pmatrix}p(v_r)& 0 \\ \star & M_2(p')  \end{pmatrix}\begin{pmatrix} 0 \\ m' \end{pmatrix}=0, \]
contradicting the infinitesimal rigidity of $(G/H,p^*)$.

(b) A similar proof holds. We choose a generic $uv$-coincident framework $(G,p)$, a vector $m\in \ker R_\Y(G,p)$ which is not an infinitesimal isometry of $\mathbb{R}^3$, and $uv$-coincident realisations $(G,p')$ and $(G/H,p^*)$ as above. We then use the facts that $H$ is rigid on $\Y$ and $G/H$ is $uv$-rigid on $\Y$ to obtain a contradiction.
\eop

We next consider the 0-$uv$-reduction, 1-$uv$-reduction, $uv$-$K_4$-contraction and $uv$-4-cycle contraction operations.

\begin{lemma} \label{lem:deg3}
Let \(G=(V,E)\) be a graph and let $u,v\in V$ be distinct vertices.
Suppose that $|E|=2|V|-2$, $G$ is independent in
\(\mathcal{M}_{uv}\), and \(d(w)\geq 3\) for all $w\in V$. Then
either there is a vertex \(z\in V-\{u,v\}\) with \(d(z)=3\) and
\(|N(z)\cap \{u,v\}|\leq 1\) or there is a $4$-cycle in $G$ which
contains both $u$ and $v$.
\end{lemma}

\proof Since $|E|=2|V|-2$ and $d(w)\geq 3$ for all $w\in V$, there
are at least 4 vertices of degree 3. Since $G$ is independent in
$\mathcal{M}_{uv}$, $G$ has at most two vertices which are adjacent
to both $u$ and $v$. Hence, if
there is no vertex $z\in V-\{u,v\}$ with \(d(z)=3\) and
\(|N(z)\cap \{u,v\}|\leq 1\), then the vertices of degree $3$ must
induce a $C_4$ in $G$ which contains both $u$ and $v$.
\eop

\begin{lemma}\label{lem:uvc4}
Let \(G=(V,E)\) be a graph and let $u,v\in V$ be distinct vertices.
Suppose that $G$ is independent in  \(\mathcal{M}_{uv}\), and there
are vertices \(a,b\) such that \(a,u,b,v\) is a cycle in \(G\). Then
the $uv$-4-cycle contraction which merges $u$ and $v$ results in a
simple graph $G'$ which is $(2,2)$-sparse.
\end{lemma}

\proof The independence of \(G\) in $\M_{uv}$ implies that there is
no vertex other than \(a,b\) that is adjacent with both \(u\) and
\(v\). Thus \(G'\) is simple.  Suppose \(G\) is not $(2,2)$-sparse.
Then there exists a $(2,2)$-tight set \(X\) in \(G\) that contains
\(u,v\) and exactly one of \(a\) and \(b\), say $a$.
Let \(\{X,\{u,v,b\}\}=\mathcal{H}\). Then
\(i(\mathcal{H})=2|X|-2+2\) and \(\val(\mathcal{H})=2|X|-2+3-2\)
which contradicts the independence of \(G\) in \(\mathcal{M}_{uv}\).
\eop

\begin{lemma} \label{lem:red}
Let $G=(V,E)$ be a graph and let $u,v\in V$ be distinct vertices.
Suppose that \(G\) is independent in \(\mathcal{M}_{uv}\) and let
\(z\in V-\{u,v\}\) with $N(z)=\{v_1,v_2,v_3\}$ and $|N(z)\cap
\{u,v\}|\leq 1$. Then either:\\
(a) there is a 1-reduction at $z$ which leads to
a graph which is independent in $\mathcal{M}_{uv}$, or\\
(b) $z$ and its neighbours induce a copy of $K_4$ in $G$, or\\
(c)  $v_i\in \{u,v\}$ and $v_jv_k\in E$ for some
$\{i,j,k\}=\{1,2,3\}$, and there is a tight $uv$-compatible family
$\{X_1,X_2,\ldots,X_k\}$ in $G$ such that $X_1=N(z)\cup \{u,v,z\}$
and $i(X_1)\geq 2|X_1|-4$.
\end{lemma}

\proof Suppose (a) does not occur. Then, for all $1\leq i<j\leq 3$,
either $v_iv_j\in E$, or
there exists a tight $uv$-compatible family $\mathcal{H}_{ij}$ in $G-z$ with $v_iv_j\in \cov(\mathcal{H}_{ij})$
or
there exists a tight set $X_{ij}$ in $G-z$ with $\{v_i,v_j\}\subset
X_{ij}$ and $\{u,v\}\not\subset X_{ij}$. If the second alternative
occurs we may assume that $\mathcal{H}_{ij}$ has been chosen to be
the unique  tight $uv$-compatible family in $G-z$ with maximal
cover. If $G[v_1,v_2,v_3]\cong K_3$ then (b) occurs. So we may
assume that $v_1v_2\notin E$.

We first show that $v_iv_j\notin E$ and that $\mathcal{H}_{ij}$
exists for some $1\leq i < j \leq 3$. Suppose $\mathcal{H}_{12}$
does not exist. Then $X_{12}$ exists. If $v_3\in X_{12}$ then
$X_{12}+ z$ contradicts the independence of $G$ in
$\mathcal{M}_{uv}$. Hence $v_3\notin X_{12}$. If $v_1v_3,v_2v_3\in
E$ then $X_{12}\cup \{v_3,z\}$ contradicts the independence of $G$
in $\mathcal{M}_{uv}$. Hence suppose that $v_1v_3\notin E$. If
$X_{13}$ exists, then $X_{12}\cup X_{13}\cup \{z\}$ contradicts the
independence of $G$ in $\mathcal{M}_{uv}$. Hence $\mathcal{H}_{13}$
exists. Relabeling if necessary we assume that
$\mathcal{H}_{12}=\{X_1,X_2,\dots,X_k\}$ exists.

Since $v_1v_2\in \cov(\mathcal{H}_{12})$ we have $v_1,v_2\in X_i$
for some $1\leq i\leq k$.
If $v_3\in X_i$ then $|X_i|\geq 4$, since $|N(z)\cap \{u,v\}|\leq
1$, and the $uv$-compatible family obtained from $\cH_{12}$ by
replacing $X_i$ by $X_i+z$ will contradict the independence of $G$
in $\M_{uv}$. Hence $v_3\not\in X_i$.

Suppose that $\{v_1,v_2\}\cap \{u,v\}=\emptyset$. Then $|X_i|\geq 4$
and neither $v_1v_3$ nor $v_2v_3$ are covered by $\cH_{12}$. The
maximality of $\cov(\cH_{12})$ now implies that $\cH_{13}$ and
$\cH_{23}$ do not exist.
If $v_1v_3,v_2v_3\in E$, then the $uv$-compatible family obtained
from $\cH_{12}$ by replacing $X_i$ by $X_i+v_3$ will contradict the
maximality of $\cov(\mathcal{H}_{12})$.
Relabeling if necessary, we may suppose that $v_1v_3\notin E$, and
hence $X_{13}$ exists. Then $X_i\cap X_{13}\neq \emptyset$,
$|X_i|\geq 4$, $|X_{13}|\geq 4$ and $v_1v_3\in \cov(X_{13})\sm
\cov(\cH_{12})$. This contradicts the final part of Lemma
\ref{tightsystem}. Hence $\{v_1,v_2\}\cap \{u,v\}\neq\emptyset$ and
we may assume, without loss of generality, that $u=v_1$.

If  $v_3\not\in V(\cH_{12})$, then Lemma \ref{tightsystem}(c)
implies that $v_1v_3\not\in E$ and hence $X_{13}$ exists. This
contradicts the final part of Lemma \ref{tightsystem} since $u\in
X_{13}\cap X_i$. Hence $v_3\in X_j$ for some $X_j\in \cH_{12}-X_i$.
The final part of Lemma \ref{tightsystem} now implies that $X_{23}$
does not exist and hence $v_2v_3\in E$.

Let $X=X_i\cup X_j\cup \{z\}$ and
$\cH=(\cH_{12}\sm\{X_i,X_j\})\cup\{X\}$. Then the facts that $G$ is
independent in $\M_{uv}$ and $\cH_{12}$ is tight imply that
$|X_i|=3=|X_j|$ (so $X=N(z)\cup \{u,v,z\}$), and that $\cH$ is a
tight $uv$-compatible family in $G$ with $i(X)\geq 2|X|-4$. \eop

\begin{lemma} \label{lem:mintight}
Let $G=(V,E)$ be a graph and let $u,v\in V$ be distinct vertices.
Suppose that \(G\) is independent in \(\mathcal{M}_{uv}\),
$\cH=\{X_1,X_2,\ldots,X_k\}$ is a tight $uv$-compatible family in
$G$ and that
$\cH-X_i$ is not tight for all $1\leq i\leq k$. Then either:\\
(a) $k=1$ and $X_1$ is tight; \\
(b) $k=2$, $|X_1|=|X_2|=3$ and $i(X_1)=i(X_2)=2$; \\
(c) $k=2$, $|X_1|\geq 4$, $i(X_1)=2|X_1|-3$, $|X_2|=3$ and $i(X_2)=2$; or\\
(d) $k=2$, $|X_i|\geq 4$ and $i(X_i)=2|X_i|-3$ for all $i\in \{1,2\}$.
\end{lemma}
\proof We have $i(\cH-X_i)<\val(\cH-X_i)$, and hence $i(X_i)\geq
2|X_i|-3$ if $|X_i|\geq 4$ and $i(X_i)= 2$ if $|X_i|= 3$. The fact
that $G$ is independent in $\mathcal{M}_{uv}$ and $\cH$ is tight now imply that $k=1$ or
$2$ and that the sets in $\cH$ satisfy the assertions in the lemma.
\eop

Note that if alternative (d) holds then $X_1\cup X_2$ is tight so we can reduce
to alternative (a).

\begin{lemma} \label{lem:k4_1}
Let $G=(V,E)$ be a graph and let $u,v\in V$ be distinct vertices.
Suppose that \(G\) is independent in \(\mathcal{M}_{uv}\) and that
there exists a subgraph $H$ of $G$ isomorphic to $K_4$.
Then either:\\
(a) there is a vertex $x\in V-V(H)$ such that $|N(x)\cap V(H)|= 2$,\\
(b) $|V(H)\cap \{u,v\}|=1=|N(V(H))\cap \{u,v\}|$,\\
(c) there is a tight $uv$-compatible family $\{X_1,X_2,\ldots,X_k\}$ in $G$ such
that $X_1=V(H)\cup\{u,v\}$, $|X_1|=6$ and $i(X_1)=8$, \\
(d) there is a tight $uv$-compatible family $\{X_1,X_2,\ldots,X_k\}$ in $G$ such
that $X_1=V(H)\cup\{u,v,a\}$ for some $a\in V-(V(H)\cup \{u,v\})$, $|X_1|=6$ and $i(X_1)=8$, or\\
(e) the  contraction of $H$ gives a graph
$G'$ which is independent in \(\mathcal{M}_{uv}\).
\end{lemma}

\proof
Let $G'$ be the result of a $K_4$ contraction applied to $H$ with \(w\) being the contracted vertex. Suppose that (a), (b) and (e) fail.
It is easy to check that $G'$ is $(2,2)$-sparse. Since (a) fails, $G'$ is simple. Since (b) fails, $uv\notin E(G')$. Since (e), there is a $uv$-compatible family $\mathcal{H}=\{X_1,X_2,\ldots,X_k\}$ for which \(\val(\mathcal{H})<i(\mathcal{H})\) and \(w\in V(\mathcal{H})\). Without loss of generality we may assume \(w\in X_1\). If \(|X_1|\geq4\) then we get a contradiction as the $uv$-compatible family \(\mathcal{H}'=\{(H_1-w)\cup V(H),X_2,\dots,X_k\}\) of \(G\) violates independence. If \(|X_1|=3\) and $V(H)\cap \{u,v\}=\emptyset$ then \(\mathcal{H}'\) is the $uv$-compatible family described in (c). Finally if \(|X_1|=3\) and $V(K_4)\cap \{u,v\}=u$ then \(\mathcal{H}''=\{V(H)\cup\{u,v,a\},X_2,\dots,X_k\}\) is the $uv$-compatible family described in (d).
\eop

\begin{lemma} \label{lem:4-cycle}
Let $G=(V,E)$ be a graph and let $u,v\in V$ be distinct vertices.
Suppose that \(G\) is independent in \(\mathcal{M}_{uv}\), $z\in
V-\{u,v\}$ is a vertex of degree 3 with $N(z)=\{v_1,v_2,v_3\}$,
$|N(z) \cap \{u,v\}|\leq 1$ and $G[N(z)+z]$ is isomorphic to $K_4$.
Suppose further that there is a vertex $x\in V-\{z,v_1,v_2,v_3\}$
such that $N(x)\cap N(z)= \{v_2,v_3\}$ and $\{v_1,x\}\neq \{u,v\}$.
Then the $uv$-4-cycle contraction operation which contracts $x$ and $z$
into a single vertex $x$ leads to a graph $G'$ which is independent
in \(\mathcal{M}_{uv}\).
\end{lemma}

\proof Suppose $G'$ is not independent in $\mathcal{M}_{uv}$. Since
$G'=G-z+v_1x$ and $xv_1\notin E$, there exists either a tight
$uv$-compatible family $\mathcal{H}$ in $G-z$ with $xv_1\in
\cov(\mathcal{H})$, or a tight set $X$ in $G-z$ with
$\{x,v_1\}\subset X$. Set
$Y=\{z,v_1,v_2,v_3,x\}$. Then $Y$ is tight in $G$.

Suppose $X$ exists. Then $X\cup Y$ and $X\cap Y$ are tight by Lemma
\ref{lem:sets}. Since $\{v_1,x\} \subseteq X\cap Y$ and no
proper subset of $Y$ containing $v_1$ and $x$ is tight, we have
$X\cap Y=Y$. This implies that $z\in X$ contradicting the choice of
$X$. Hence $\mathcal{H}=\{X_1,X_2,\dots,X_k\}$ exists.

Since $xv_1\in \cov(\mathcal{H})$, we may assume, without loss of
generality, that $x,v_1\in X_1$. Then $x,v_1\in X_1\cap Y$. Since \(|\{u,v\}\cap Y|\leq 1\) by the hypotheses of the lemma, Lemma \ref{lem:set2intersection} implies that $X_1\cap Y$ is tight. Since no
proper subset of $Y$ containing $v_1$ and $x$ is tight we have
$X_1\cap Y=Y$. This implies that $z\in X_1$ and contradicts the choice of $\cH$.
\eop

\begin{theorem}
\label{thm:char}
Let \(G=(V,E)\) be a graph and let $u,v\in V$ be distinct vertices. Then \(G\) is independent in $\mathcal{R}_{uv}^\Y$ if and only if \(G\) is independent in \(\mathcal{M}_{uv}\).
\end{theorem}

\proof
Necessity follows from Lemma \ref{lem:necessary}. Now suppose that \(G\) is independent in \(\mathcal{M}_{uv}\). We prove that \(G\) is independent in $\mathcal{R}_{uv}^\Y$ by induction on $|V|$. It is straightforward to check that \(G\) is independent in $\mathcal{R}_{uv}^\Y$ when $|V| \leq 4$.  Hence we may assume that $|V| \geq 5$. By extending $|E|$ to a base of $\mathcal{M}_{uv}(K_{|V|})$ if necessary, we may also assume that $|E|=2|V|-2$.

\smallskip\noindent
\textbf{Case 1. {\boldmath $G$ contains a vertex of degree
2.}} First suppose that $u$ has degree 2. Then $G-u$ is
$(2,2)$-sparse. Hence, by Theorem \ref{thm:cylinderlaman},
$R_\Y(G-u,p)$ has linearly independent rows for any generic $p$. We can now use Lemma \ref{lem:0ext} to show that $G$ is independent in $\mathcal{R}_{uv}^\Y$.

Now, suppose that there is a vertex \(w\in V-\{u,v\}\) with \(d(w)=2\). Let \(N(w)=\{a,b\}\). Clearly, \(a\neq b\) holds. If \(\{a,b\}=\{u,v\}\) then let \(\mathcal{H}=\{\{u,v,w\},\{V-w\}\}\), where $|V-w|\geq 4$. We have
\[2|V|-2=|E|=i_{E}(\mathcal{H})\leq \val(\mathcal{H})=2\cdot3-3+2(|V|-1)-2-2=2|V|-3,\]
a contradiction. Hence \(\{a,b\}\neq\{u,v\}\), which implies that the 0-\(uv\)-reduction operation can be applied at \(w\) to obtain a graph \(G'=(V-w,E')\) that is independent in \(\mathcal{M}_{uv}\) and satisfies \(|E'|=2|V-w|-2\). By induction, \(G'\)  is independent in $\mathcal{R}_{uv}^\Y$. Now Lemma \ref{lem:ext} implies that \(G\) is independent in $\mathcal{R}_{uv}^\Y$.

\smallskip\noindent
\textbf{Case 2. {\boldmath There is a 4-cycle in $G$ containing $u$
and $v$.}} By Lemma \ref{lem:uvc4}, we may apply a
$uv$-4-cycle-contraction (contracting $u$ and $v$) to obtain a graph
$H$ which is simple and $(2,2)$-sparse. Theorem
\ref{thm:cylinderlaman} implies that any generic realisation $(H,p)$
on $\Y$ is infinitesimally rigid. Now we can use Lemma
\ref{lem:4cycleind} to show that \(G\) is independent in
$\mathcal{R}_{uv}^\Y$.

\smallskip

Henceforth we assume that Cases 1 and 2 do not occur.

\smallskip\noindent
\textbf{Case 3. {\boldmath There is a proper tight set $X$
containing $u$ and $v$.}} Since Case 1  does not occur, we may
suppose $X$ is a maximal proper tight set (where proper means $X\neq
V$ and maximal means there is no vertex $w\in V-X$ with more than
one neighbour in $X$). Now by the maximality of $X$, $G/X$ is simple
and $|V-X|\geq 3$. Hence $G/X$ is $(2,2)$-tight. Theorem
\ref{thm:cylinderlaman} implies that any generic framework
$(G/X_1,p)$ on $\Y$ is infinitesimally rigid. We may now apply Lemma
\ref{lem:extension}(a) to show that \(G\) is independent in
$\mathcal{R}_{uv}^\Y$.

\smallskip

Henceforth we may assume that Case 3 does not occur.

\smallskip\noindent
\textbf{Case 4. {\boldmath There is a degree three vertex $z$ in $G$ which is contained in a subgraph $H\cong K_4$,  and a vertex $x\in V-V(H)$ such that $|V(H)\cap N(x)|=2$.}} If $\{u,v\}\not\subset V(H)\cup \{x\}$ then we may apply Lemma \ref{lem:4-cycle} to find a graph $G'$ which is independent in $\mathcal{M}_{uv}$. We can now use Corollary \ref{cor:4cycleuv} to show that \(G\) is independent in $\mathcal{R}_{uv}^\Y$. Thus we may suppose that $\{u,v\} \subset V(H)\cup \{x\}$. Then $H\cup \{x\}$ is tight. This contradicts the assumption that Case 1 (if $H\cup \{x\}=V$) or Case 3 (if $H\cup\{x\}\neq V$) do not occur.

\smallskip

A vertex $z$ of degree 3 in $G$ is \emph{bad} if either
\begin{itemize}
\item
$z\in \{u,v\}$, or
\item $z$ is adjacent to both $u$ and $v$,or
\item
$X=N(z)\cup \{u,v,z\}$ satisfies alternative (c) of Lemma
\ref{lem:red} and $i(X)\geq 2|X|-3$, or
\item
$z$ belongs to a subgraph $H\cong K_4$ satisfying alternative (b) of
Lemma \ref{lem:k4_1}.
\end{itemize}
Otherwise we say that $z$ is \emph{good}.

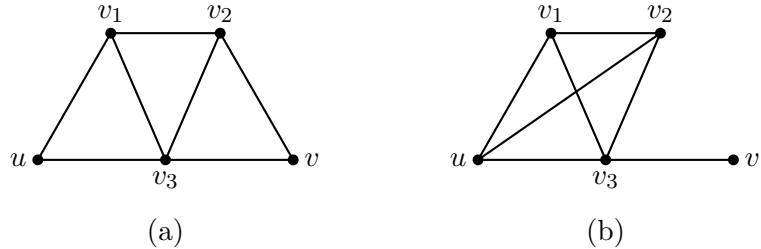
\begin{figure}[htp]
\begin{center}
     \begin{tikzpicture}[very thick,scale=.48]
\filldraw (0,0) circle (3pt)node[anchor=east]{$u$};
\filldraw (3.5,0) circle (3pt)node[anchor=north]{$v_3$};
\filldraw (2,3.5) circle (3pt)node[anchor=south]{$v_1$};
\filldraw (5,3.5) circle (3pt)node[anchor=south]{$v_2$};
\filldraw (7,0) circle (3pt)node[anchor=west]{$v$};

 \draw[black,thick]
(5,3.5) -- (3.5,0) -- (2,3.5) -- (0,0) -- (3.5,0) -- (7,0) -- (5,3.5) -- (2,3.5);

\node [rectangle, draw=white, fill=white] (b) at (3.5,-2) {(a)};
\end{tikzpicture}
       \hspace{1.2cm}
    \begin{tikzpicture}[very thick,scale=.48]
\filldraw (0,0) circle (3pt)node[anchor=east]{$u$};
\filldraw (3.5,0) circle (3pt)node[anchor=north]{$v_3$};
\filldraw (2,3.5) circle (3pt)node[anchor=south]{$v_1$};
\filldraw (5,3.5) circle (3pt)node[anchor=south]{$v_2$};
\filldraw (7,0) circle (3pt)node[anchor=west]{$v$};

 \draw[black,thick]
(5,3.5) -- (3.5,0) -- (2,3.5) -- (0,0) -- (3.5,0) -- (7,0);

\draw[black,thick]
(0,0) -- (5,3.5) -- (2,3.5);

  \node [rectangle, draw=white, fill=white] (b) at (3.5,-2) {(b)};
         \end{tikzpicture}
\end{center}
\vspace{-0.3cm}
\caption{The two alternatives for $H$.}
\label{fig1}
\end{figure}

\smallskip\noindent
\textbf{Case 5. All degree three vertices are bad}. We may use Lemma \ref{lem:deg3} and the fact that Case 2 does not occur to deduce there exists a degree three vertex $v_1\in V\sm \{u,v\}$ with $|N(v_1)\cap \{u,v\}|\leq 1$. Since $v_1$ is bad $X=N(v_1)\cup\{u,v,v_1\}$ satisfies alternative (c) of Lemma \ref{lem:red} and $i(X)\geq 2|X|-3$, or $v_1$ belongs to a subgraph $H\cong K_4$ satisfying alternative (b) of Lemma \ref{lem:k4_1}. If the first alternative occurs then we may use the facts that $G$ is independent in $\M_{uv}$ and Case 2 does not occur to deduce that $i(X)=2|X|-3=7$. It follows that, in both cases, we may relabel the vertices of $H=G[N(v_1)\cup \{u,v,v_1\}]$ such that $H$ is one of the graphs shown in Figure \ref{fig1}.

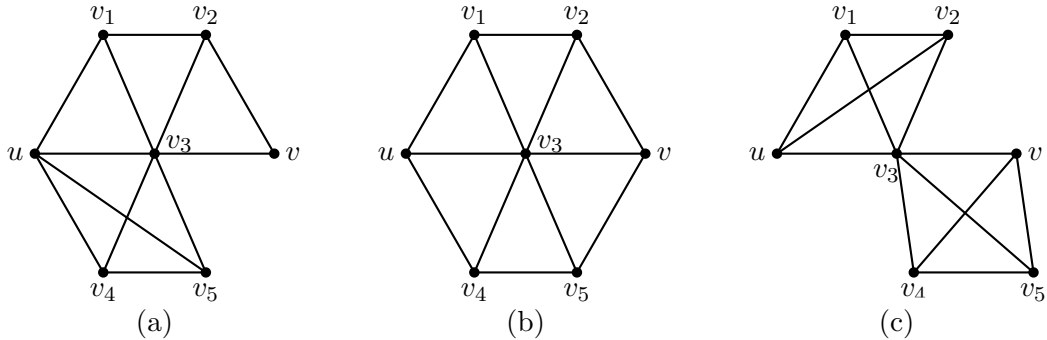
\begin{figure}[htp]
\begin{center}
   \begin{tikzpicture}[very thick,scale=.45]
\filldraw (0,0) circle (3pt)node[anchor=east]{$u$};
\filldraw (3.5,0) circle (3pt);
\filldraw (2,3.5) circle (3pt)node[anchor=south]{$v_1$};
\filldraw (5,3.5) circle (3pt)node[anchor=south]{$v_2$};
\filldraw (7,0) circle (3pt)node[anchor=west]{$v$};
\filldraw (2,-3.5) circle (3pt)node[anchor=north]{$v_4$};
\filldraw (5,-3.5) circle (3pt)node[anchor=north]{$v_5$};

\filldraw (3.5,.3) circle (0pt)node[anchor=west]{$v_3$};

 \draw[black,thick]
(5,3.5) -- (3.5,0) -- (2,3.5) -- (0,0) -- (3.5,0) -- (7,0) -- (5,3.5) -- (2,3.5);

\draw[black,thick]
(5,-3.5) -- (0,0) -- (2,-3.5) -- (3.5,0) -- (5,-3.5);

\draw[black,thick]
(2,-3.5) -- (5,-3.5);

\node [rectangle, draw=white, fill=white] (b) at (3.5,-5) {(a)};
\end{tikzpicture}
       \hspace{.5cm}
     \begin{tikzpicture}[very thick,scale=.45]
\filldraw (0,0) circle (3pt)node[anchor=east]{$u$};
\filldraw (3.5,0) circle (3pt);
\filldraw (2,3.5) circle (3pt)node[anchor=south]{$v_1$};
\filldraw (5,3.5) circle (3pt)node[anchor=south]{$v_2$};
\filldraw (7,0) circle (3pt)node[anchor=west]{$v$};
\filldraw (2,-3.5) circle (3pt)node[anchor=north]{$v_4$};
\filldraw (5,-3.5) circle (3pt)node[anchor=north]{$v_5$};

\filldraw (3.5,.3) circle (0pt)node[anchor=west]{$v_3$};

 \draw[black,thick]
(5,3.5) -- (3.5,0) -- (2,3.5) -- (0,0) -- (3.5,0) -- (7,0) -- (5,3.5) -- (2,3.5);

\draw[black,thick]
(0,0) -- (2,-3.5) -- (3.5,0) -- (5,-3.5) -- (7,0);

\draw[black,thick]
(2,-3.5) -- (5,-3.5);

\node [rectangle, draw=white, fill=white] (b) at (3.5,-5) {(b)};
\end{tikzpicture}
       \hspace{.5cm}
    \begin{tikzpicture}[very thick,scale=.45]
\filldraw (0,0) circle (3pt)node[anchor=east]{$u$};
\filldraw (3.5,0) circle (3pt);
\filldraw (2,3.5) circle (3pt)node[anchor=south]{$v_1$};
\filldraw (5,3.5) circle (3pt)node[anchor=south]{$v_2$};
\filldraw (7,0) circle (3pt)node[anchor=west]{$v$};
\filldraw (4,-3.5) circle (3pt)node[anchor=north]{$v_4$};
\filldraw (7.5,-3.5) circle (3pt)node[anchor=north]{$v_5$};

\filldraw (3.2,0) circle (0pt)node[anchor=north]{$v_3$};

 \draw[black,thick]
(5,3.5) -- (3.5,0) -- (2,3.5) -- (0,0) -- (3.5,0) -- (7,0) -- (7.5,-3.5) -- (4,-3.5) -- (3.5,0) -- (7.5,-3.5);

\draw[black,thick]
(0,0) -- (5,3.5) -- (2,3.5);

\draw[black,thick]
(7,0) -- (4,-3.5);

  \node [rectangle, draw=white, fill=white] (b) at (3.5,-5) {(c)};
         \end{tikzpicture}
\end{center}
\vspace{-0.3cm}
\caption{The three alternatives for $G$.}
\label{fig2}
\end{figure}

The fact that $G$ is $(2,2)$-sparse implies that, in both cases, there exists a (necessarily bad) degree three vertex $v_4\in V\sm V(H)$. Since Case 2 does not occur, $v_4$ is not adjacent to both $u$ and $v$. Hence $v_4$  also belongs to a subgraph $H'$ which is isomorphic to one of the graphs shown in Figure \ref{fig1}. Since Case 2 does not occur, $v_4\in V(H')$. Since  $G$ is $(2,2)$-sparse, $V(H)\cap V(H')=\{u,v,v_3\}$. Now $H\cup H'$ is one of the graphs shown in
Figure \ref{fig2}. Since all three graphs are tight, we may use the fact that Case 3 does not occur to deduce that $G=H \cup H'$. The fact that Case 1 does not occur now tells us that $G$ is not the graph in Figure \ref{fig2}(a). The graph in Figure \ref{fig2}(b) cannot be equal to $G$ since $X=N(v_1)\cup \{u,v,v_1\}$ does not belong to a tight $uv$-compatible family (so $v_1$ is not bad). Hence $G$ is as shown in Figure \ref{fig2}(c).

We will complete the discussion of this case by showing that $G$ is minimally $uv$-rigid on $\mathcal{Y}$. Let $(G, p)$ be a generic $uv$-coincident realisation of $G$ on $\Y$ and $m$ be an infinitesimal motion of $(G, p)$ with $m (u)=0$. Since $K_4$ is rigid, $m (w)=0$ for all $w\in V (H)-v$. In particular $m (v_3)=0$
and hence $m (w)=0$ for all $w\in V$.

\smallskip\noindent
\textbf{Case 6. None of the previous cases occur}. Let $z_1,z_2,\ldots,z_k$ be the good degree three vertices in $G$. If the edge set of some $1$-reduction of $G$ at $z_i$ is independent in $\M_{uv}$ then we may apply induction to the reduced graph and then apply Lemma \ref{lem:ext} to deduce that $G$ is independent in
$\R^\Y_{uv}$. Hence we may assume  that alternative (b) or (c) of Lemma \ref{lem:red} holds for $z_i$. Similarly, if alternative (b) of Lemma \ref{lem:red} holds for $z_i$ and $z_i$ is contained in a $K_4$-subgraph whose contraction results in a graph which is independent in $\M_{uv}$, then we may apply induction to the reduced graph and then apply Lemma \ref{lem:extension}(b) to deduce that $G$ is independent in $\R^\Y_{uv}$.  Since Case 4 does not occur, it follows that, for every good degree 3 vertex $z_i$, there exists a tight $uv$-compatible family $\cH_i$ as described in alternative (c) or (d) of Lemma \ref{lem:k4_1}, or alternative (c) of Lemma \ref{lem:red}. We may assume that the first alternatives holds for all $1\leq i\leq k$ and that the second alternative holds for
$l+1\leq i\leq k$. Choose a subgraph \(H_i\cong K_4\) that contains $z_i$ and satisfies alternative (c) or (d) of Lemma \ref{lem:k4_1} for each $1\leq i\leq k$ and let  $X_i$ be the element of $\cH_i$ which contains $V(H_i)$. Let  \(X_i=\{z_i,u,v\}\cup N(z_i)\) for each $l+1\leq i \leq k$. With these definitions we have \(i(X_i)=2|X_i|-4\) for all \(1\leq i\leq k\).

Let \(X=\bigcup_{i=1}^k X_i\). We will show by induction that \(i(X)\geq2|X|-4\). Suppose that we have \(i(X')\geq2|X'|-4\) for some $X'=\bigcup_{i=1}^s X_i$ and some $1\leq s \leq k$. If \(i(X'\cup X_{s+1})\leq 2|X'\cup X_{s+1}|-5\), then Lemma \ref{lem:2sets}(a) implies that \(i(X'\cap X_{s+1})\geq 2|X'\cap X_{s+1}|-3\), contradicting the fact that no subset of \(X_{s+1}\) that contains \(u,v\) (in each of the three possibilities for $X_{s+1}$) satisfies this inequality.

We may apply Lemma \ref{lem:mintight} to a minimal tight $uv$-compatible subfamily of $\mathcal{H}_i$ for all $1\leq i\leq k$,
and use the facts that  Cases 2 and 3 do not occur to deduce that alternative (c) of Lemma \ref{lem:mintight} must hold. Hence there exist sets \(Y_i\) and $\{u,v,y_i\}$ in $\mathcal{H}_i- X_i$ with \(i(Y_i)=2|Y_i|-3\) and $i(\{u,v,y_i\})=2$. Lemma \ref{lem:intersect3}(b) implies that \(Y_i\cap X_i=\{u,v\}=Y_i\cap\{u,v,y_i\}\) for all \(1\leq i\leq k\).
The fact that we are not in Case 2 also implies that $y_i=y_j=y$, say, for all $1\leq i \leq j \leq k$.
Let $Y=\cap_{i=1}^k Y_i$.
Then $Y\cap X=\{u,v\}$ and $y\not\in Y$. We can now use Lemma \ref{lem:2sets}(a) and the fact that Case 3 does not occur to prove inductively that \(i(Y)=2|Y|-3\).

Let \(W=V\setminus X\). Since \(i(W)\leq2|W|-2\) there is an integer
\(t\) for which \(i(W)=2|W|-2-t\). Since $i(Y)=2|Y|-3$ and $G$ is
$(2,2)$-sparse, there are at least 3 edges from $Y\sm \{u,v\}$ to
$\{u,v\}$. Since \(Y\sm\{u,v\}\subseteq W\), $y\in W\sm Y$  and
there are two edges from $y$ to $\{u,v\}$, we have at least five
edges between \(\{u,v\}\) and \(W\). Note that the definition of $X$
tells us that all degree 3 vertices in $W$ are bad.

Suppose that no (bad) degree three vertex $z\in W$ is contained in a
set $X\subseteq V$ which satisfies alternative (c) of Lemma
\ref{lem:red} and has $i(X)\geq 2|X|-3$, or a subgraph \(H\cong
K_4\) that satisfies alternative (b) of Lemma \ref{lem:k4_1}. Then
every (bad) degree three vertex in $W$ is adjacent to both $u$ and
$v$. Since Case 2 does not occur we have at most one degree three
vertex in $W$. Since \(i(Y)\geq 2|Y|-4\), we have
\(|E|-|E(Y)|-|E(W)|\leq 4+t\). The total degree of the vertices in
\(W\) is at most \(2(2|W|-2-t)+4+t=4|W|-t\). Since there is at most
one degree three vertex in \(W\), \(t\leq 1\). If $t=0$, then $W$ is
tight and \(W+u+v\) violates sparsity. Hence $t=1$ and
\(W+u+v\) is a proper tight set which contradicts the fact that Case 3 does not occur.

Now consider the case when there is a (bad) degree three vertex
$z\in W$ which is contained in a set $X\subseteq V$ which satisfies
alternative (c) of Lemma \ref{lem:red} and has $i(X)\geq 2|X|-3$, or
a subgraph \(H_1\cong K_4\) that satisfies alternative (b) of Lemma
\ref{lem:k4_1}.
Then $H=G[N(z)\cup \{u,v,z\}]$ is isomorphic to one of the graphs
shown in Figure \ref{fig1}, with $v_1=z$. Since Case 2 does not
occur we have $v_3=y$. The facts that $y\not\in Y$ and no
$Z\subseteq V(H)-z$, with $\{u,v\}\subset Z$, has $i(Z)=2|Z|-3$ imply that  $Y\cap V(H)$ is a
proper subset of both $Y$ and $V(H)$.
Lemma \ref{lem:2sets}(a)
now implies that \(Y\cup V(H)\) is tight. Since \(Y\cup V(H)\neq V\)
this contradicts the fact that Case 3 does not occur. \eop

\subsection{A deletion-contraction characterisation of $uv$-rigidity}

\begin{theorem} \label{merevsegkov}
Let \(G=(V,E)\) be a graph and let \(u,v\in V\) be distinct vertices. Then \(G\) is \(uv\)-rigid if and only if \(G-uv\) and \(G/uv\) are both rigid.
\end{theorem}

\proof
Necessity follows from the fact that an infinitesimally rigid $uv$-coincident realisation of $G$ is an infinitesimally rigid realisation of $G-uv$, and also gives rise to an infinitesimally rigid realisation of $G/uv$ by (\ref{eq:kernel}).

To prove sufficiency, suppose, for a contradiction, that \(G-uv\) and \(G/uv\) are both rigid but \(G\) is not \(uv\)-rigid. By Theorems \ref{matroid} and \ref{thm:char} this implies that there is a thin cover $\mathcal{K}$ of $G-uv$ with \(\val(\mathcal{K})\leq 2|V|-3\). If $\mathcal{K}$ consists of subsets of $V$
only, then $r(G-uv)\leq 2|V|-3$ follows, which contradicts the fact that $G-uv$ is rigid.

Hence \(\mathcal{K}=\{\mathcal{H},H_{1},\dots,H_{k}\}\),
where \(\mathcal{H}=\{X_{1},\dots,X_{l}\}\) is a \(uv\)-compatible family.
Contract the vertex pair $u,v$ in $G$ into a new vertex $z_{uv}$.
This gives rise to a cover
$$\mathcal{K}'=\{X_{1}',\dots, X_{l}',H_{1},\dots,H_{k}\}$$
of $G/uv$, where $X_j'$ is obtained from $X_j$ by replacing
$u,v$ by $z_{uv}$, for $1\leq j\leq l$.
Then we obtain
\[\sum_{i=1}^{k}(2|H_{i}|-t_{H_i})+\sum_{j=1}^{l}(2|X_{j}'|-t(X_j'))\leq
\sum_{i=1}^{k}(2|H_{i}|-t_{H_i})+\]
\[+\sum_{j=1}^{l}(2|X_{j}|-t(X_j))-2l=
\val(\mathcal{K})-2\leq 2|V|-3-2=2(|V|-1)-3,\]
which implies  that \(G/uv\) is not rigid, a contradiction. This
completes the proof.
\eop

A similar proof can be used to verify the following more general result:

\begin{theorem}
\label{rankformula}
Let \(G=(V,E)\) be a graph and let \(u,v\in V\) be distinct vertices. Then \(r_{uv}(G)=\min\{r(G-uv),r(G/uv)+2\}.\)
\end{theorem}

Theorems \ref{merevsegkov} and \ref{rankformula} show that the polynomial-time algorithms for computing the rank of a count matroid (see e.g. \cite{BJ,LS}) can be used to test whether $G$ is $uv$-rigid, or more generally, to compute $r_{uv}(G)$.

\section{Vertex splitting and global rigidity}
\label{sec:global}

Suppose $G=(V,E)$ is a graph with $V=\{v_1,v_2.\ldots,v_n\}$ and $(G,p)$ is a realisation of $G$ on a family of (not necessarily distinct) concentric cylinders $\mathcal{Y}=\Y_1\cup \Y_2\cup\ldots\cup \Y_n$ such that $p(v_i)\in \Y_i$ for $1\leq i \leq n$. We say that $(G,p)$ is \emph{globally rigid} if every equivalent framework $(G,q)$ on $\mathcal{Y}$, with $q(v_i)\in \Y_i$ for all $1\leq i \leq  n$,  is congruent to $(G,p)$.

Let $G=(V,E)$ be a graph and $v_1$ be a vertex of $G$ with neighbours $v_2,v_3,\ldots,v_{t}$. A \emph{vertex split} of $G$ at $v_1$ is a graph $\hat G$ which is obtained from $G$ by deleting the edges $v_1v_2,v_1v_3,\ldots,v_1v_{k}$ and adding a new vertex $v_0$ and new edges $v_0v_1,v_0v_2,\ldots,v_0v_{k}$, for some $2 \leq k\leq t$. We will refer to the new edge $v_0v_1$ as the \emph{bridging edge} of the vertex split. We will show in this section that a vertex splitting operation, in which the bridging edge is redundant, preserves generic global rigidity on the cylinder.

Given a map $p:V\rightarrow \mathbb{R}^{3n}$, there is a unique family of concentric cylinders $\Y$ with $p(v_i)\in \Y_i$ for all $1\leq i \leq n$ as long as $p(v_i)$ does not lie on the $z$-axis for all $1\leq i \leq n$. We will refer to $\Y$ as the family of concentric cylinders induced by $p$ and denote it by $\Y^p$.
We shall need the following analogue of \cite[Theorem 13]{C&W}.

\begin{lemma}\label{lem:perturb}
If $(G,p)$ is infinitesimally rigid and globally rigid on $\Y$, then there exists an open neighbourhood $N_p$ of $p$ on $\Y$ such that for any $q\in N_p$ the framework $(G,q)$ is infinitesimally rigid and globally rigid on $\Y$.
\end{lemma}

\proof
Suppose $|V|\geq 5$ and that for any open neighbourhood $N_p$, there is a $p^*\in N_p$ such that the framework $(G,p^*)$ is not globally rigid on $\Y$. Then there is a convergent sequence $(G,p^k)$ of non-globally rigid frameworks converging to $(G,p)$.
For each framework $(G,p^k)$, let $(G,q^k)$ be an equivalent but non-congruent realisation on $\Y$. We may assume that $(G,p^k)$ and $(G,q^k)$ are in standard position (that is $p^k(v_1)=q^k(v_1)=(0,1,0)$ assuming, without loss of generality, that $r_1=1$). By the compactness of $\mathbb{R}^{3|V|}$, there is a convergent subsequence $(G,q^m)$ converging to a limiting framework $(G,q)$. As the limits of the respective sequences, $(G,q)$ must be equivalent to $(G,p)$.

If $(G,q)$ is not congruent to $(G,p)$ then we contradict the global rigidity of $(G,p)$. So $(G,p)$ and $(G,q)$ are congruent, i.e. we can transform $q$ to $p$ by a reflection in the plane $x=0$, a reflection in the plane $z=0$ or a combination of the two. We apply this same congruence to all the $(G,q^m)$ to obtain a sequence $(G,r^m)$ converging to $(G,p)$ with $(G,r^m)$ being equivalent but not congruent to $(G,p^m)$ for each $m$.

We next show that $p^m-r^m$ gives an infinitesimal motion of $(G,\frac{p^m+r^m}{2})$ on $\Y^{\frac{p^m+r^m}{2}}$. For each edge $v_iv_j$ we have
$$ \left(\frac{p^m(v_i)+r^m(v_i)}{2} - \frac{p^m(v_j)+r^m(v_j)}{2}\right) \cdot \left((p^m(v_i)-r^m(v_i))-(p^m(v_j)-r^m(v_j))\right) $$

$$ = \frac{1}{2} ((p^m(v_i)-p^m(v_j))+(r^m(v_i)-r^m(v_j)))\cdot ((p^m(v_i)-p^m(v_j))-(r^m(v_i)-r^m(v_j)))$$

$$ =\frac{1}{2} \left((p^m(v_i)-p^m(v_j))^2-(r^m(v_i)-r^m(v_j))^2\right)=0.$$
Recall that $\bar p_m(v_i)$ and $\bar r_m(v_i)$ denote the projections of $p_m(v_i)$ and $r_m(v_i)$ onto the plane $z=0$. Since $p_m(v_i)$ and $r_m(v_i)$ both lie on $\Y_i$, we have $\bar p_m(v_i)\cdot \bar p_m(v_i)=\bar r_m(v_i)\cdot \bar r_m(v_i)$. Hence for each vertex $v_i$,
$$ (\bar p_m(v_i)+\bar r_m(v_i))\cdot(\bar p_m(v_i)-\bar r_m(v_i))=0.$$

Since $p^m$ and $r^m$ are not congruent, $p^m-r^m$ is a nontrivial infinitesimal motion. This means that the rank of the rigidity matrix for each framework $(G,\frac{p^m+r^m}{2})$ is less than maximal. Since both $p^m$ and $r^m$ converge to $p$, so does $\frac{p^m+r^m}{2}$. Thus $(G,p)$ is a limit of a sequence of infinitesimally flexible frameworks and hence itself is infinitesimally flexible, a contradiction. (The fact that $(G,p)$ is infinitesimally rigid implies that the rank of $R_{\Y^q}(G,p)$ is maximum for all $q\in \mathbb{R}^{3|V|}$ sufficiently close to $p$.)
\eop

We can use this lemma and our main result to show that vertex splitting preserves global rigidity on $\Y$ under the additional assumption that the new edge is redundant.

\begin{theorem}
Let $(G,p)$ be a generic globally rigid framework on a family of concentric cylinders $\Y$. Let $\hat G$ be a vertex split of $G$ at the vertex $v_1$ with new vertex $v_0$ and suppose that $\hat G-v_0v_1$ is rigid on $\Y$. Let $\hat p(v)=p(v)$ for all $v\neq v_0$ and $\hat p(v_0)=p(v_1)$.
Then for any $q$ on $\Y$ which is sufficiently close to $\hat p$, $(\hat G,q)$ is globally rigid on $\Y$.
\end{theorem}

\proof
Since $(\hat G/v_0v_1,p)=(G,p)$ is globally rigid on $\Y$ and $p$ is generic, $\hat G/v_0v_1$ is
rigid on $\Y$. Since $G-v_0v_1$ is also rigid on $\Y$, Theorem \ref{merevsegkov} implies that $\hat G$ has a $v_0v_1$-coincident generic rigid realisation $(\hat G,\hat p)$, where $\hat p(v)=p(v)$ for all $v\neq v_0$ and $\hat p(v_0)=p(v_1)$. Since $(G,p)$ is globally rigid on $\Y$, $(\hat G,\hat p)$ is also globally rigid on $\Y$. We can now use Lemma \ref{lem:perturb} to deduce that  $(\hat G,q)$ is globally rigid on $\Y$ for all $q$ sufficiently close to $\hat p$.
\eop

\section{Concluding remarks}
\label{sec:further}

Similarly to our definition of a framework $(G,p)$ on $\Y$ we can define a framework on a family of concentric spheres $\mathcal{S}=\mathcal{S}_1\cup \mathcal{S}_2\cup \dots \cup \mathcal{S}_k$ where $\mathcal{S}_i=\{(x,y,z)\in \mathbb{R}^3:x^2+y^2+z^2=r_i\}$ and $r=(r_1,\dots,r_k)$ is a vector of positive real numbers. We can map a framework on $\mathcal{S}$ to a framework in the union of parallel (affine) planes $P_1\cup P_2\cup \dots \cup P_k$, where $P_i$ is the plane $z=r_i$ in $\mathbb{R}^3$, by central projection. In \cite{SalW,SchW} this process was shown to preserve infinitesimal rigidity for any framework on $S$. Since the projection also preserves the property that $u$ an $v$ are coincident, the problem of characterising generic rigidity for frameworks with two coincident points on concentric spheres is equivalent to the problem of characterising generic rigidity for frameworks with two coincident points on parallel planes. This latter problem can be characterised using the proof technique of Theorem \ref{thm:planeuv}. This gives us the following result.

\begin{theorem}
Let $G=(V,E)$ be a graph and let $u,v\in V$ be distinct vertices. Then $G$ is $uv$-rigid on $\mathcal{S}$ if and only if $G-uv$ and $G/uv$ are both rigid on $\mathcal{S}$.
\end{theorem}

Note that a graph $G=(V, E)$ is rigid on $\mathcal{S}$ if and only if it has rank $2|V|-3$  in the $(2, 3)$-sparse matroid by \cite[Theorem 5.1]{NOP}.

We can also replace $\Y$ with other surfaces. In particular if we choose a surface with 1 ambient rigid motion (such as the cone, hyperboloid or torus) then the analogue of Theorem \ref{thm:cylinderlaman} requires the graph to be $(2,1)$-tight \cite{NOP2}. In the $uv$-coincident case we would define the value as $\val(H)=2|H|-t_H$ where $t_H=3$ if $|H|\in \{2,3\}$ and \(H\neq\{u,v\}\), $t_H=2$ if $|H|\in\{0,4\}$ or \(H=\{u,v\}\) and $t_H=1$ if $|H|\geq 5$. We expect that, using similar techniques to Section \ref{countsec}, the appropriate count matroid can be established. However we do not know how to prove an analogue of Theorem \ref{thm:char}. To make a start on this problem would require dealing with the case when the only vertices of degree less than 4 are $u$ and $v$.


\section{Acknowledgements}
The second author was supported by the EPSRC First Grant EP/M013642/1 and by the Hungarian Scientific Research Fund (OTKA, grant number K109240).


\end{document}